\RequirePackage{fix-cm}
\documentclass[a4paper, 12pt, leqno]{article} 

\usepackage{latexsym, amsmath, amsthm, amssymb, graphicx, graphics, epsfig}
\usepackage{bm}
\usepackage{tikz-cd}

\usepackage{fullpage}

\usepackage{enumitem}

\newcommand{\comment}[1]{}

\usepackage[T1]{fontenc}
\usepackage{latexsym}
\usepackage{manfnt}
\usepackage[sans]{dsfont}
\usepackage{palatino}
\usepackage[sc, osf]{mathpazo} 
\usepackage{eulervm}
\usepackage{euscript}

\newcommand{\sub}[1]{\vspace{6pt}\textbf{#1}\hspace*{0.5em}}
\newcommand{\upskip}{\vspace{-6pt}}
\newcommand{\uppskip}{\vspace{-12pt}}

\newtheorem{thm}{Theorem}[section]
\newtheorem*{thm*}{Theorem}
\newtheorem{lem}[thm]{Lemma}
\newtheorem*{lem*}{Lemma}
\newtheorem{cor}[thm]{Corollary}
\newtheorem*{cor*}{Corollary}
\newtheorem{prop}[thm]{Proposition}
\newtheorem*{prop*}{Proposition}
\newtheorem*{claim*}{Claim}
\newtheorem*{ass*}{Assumption}

\theoremstyle{definition}

\newtheorem{rem}[thm]{Remark}
\newtheorem*{rem*}{Remark}

\newtheorem{exm}[thm]{Example}

\numberwithin{equation}{section}

\renewcommand{\proof}{\vspace{-6pt}\noindent\textit{\textbf{Proof. }}}
\newcommand{\prooff}[1]{\upskip\noindent\textit{\textbf{Proof{#1}. }}}
\renewcommand{\endproof}{$\square$}

\newcommand{\q}[1]{{``#1''}}
\newcommand{\br}[1]{\left(#1\right)}
\newcommand{\brr}[1]{\left[#1\right]}
\newcommand{\bra}[1]{\left\langle #1 \right\rangle}

\renewcommand{\leq}{\leqslant}
\renewcommand{\geq}{\geqslant}

\newcommand{\hookto}{\hookrightarrow}

\newcommand{\set}[1]{\left\{#1\right\}}
\newcommand{\sett}[2]{\left\{#1 \,|\, #2\right\}}

\newcommand{\sa}{^{\mathrm{sa}}}
\newcommand{\el}{^{\mathrm{ell}}}
\newcommand{\elsa}{^{\mathrm{ell,sa}}}

\newcommand{\eu}{_{\mathrm{eu}}}

\newcommand{\kat}{\tilde{\kappa}}

\newcommand{\Ft}{\tilde{F}}
\newcommand{\At}{\tilde{A}}

\newcommand{\Dt}{\tilde{D}}
\newcommand{\Nt}{\tilde{N}}
\newcommand{\AAt}{\tilde{\A}}

\newcommand{\xX}{x\in X}
\DeclareMathOperator{\Res}{Res}

\newcommand{\CRRsa}{\CRR\sa(H)}
\newcommand{\CRRg}{\CRR^1(H)}
\newcommand{\Bsa}{\B\sa(H)}
\newcommand{\Beusa}{\B\sa_{\mathrm{eu}}(H)}

\newcommand{\p}{\partial}

\newcommand{\pN}{\p N}

\newcommand{\Zg}{\ZZ^1}
\newcommand{\Zgg}{\bar{\ZZ}^1}
\newcommand{\Pg}{\Proj^1}

\newcommand{\Sg}{\S^1}
\newcommand{\Gg}{\Gs^1}
\newcommand{\Bg}{\B^1}
\newcommand{\Tg}{\T^1}

\newcommand{\pis}{\pi_{\S}}
\newcommand{\pig}{\pi_{\Gs}}
\DeclareMathOperator{\diag}{diag}
\newcommand{\Cl}{\mathrm{Cl}}

\newcommand{\Ho}{H^{\circ}}

\DeclareMathOperator{\one}{\mathbb{1}}
\newcommand{\chip}{\one_{[0,+\infty)}}
\newcommand{\chilp}{\one_{[\lambda,+\infty)}}

\newcommand{\chill}{\one_{(-\lambda,\,\lambda)}}
\newcommand{\chirp}{\one_{[r,+\infty)}}
\newcommand{\chirm}{\one_{(-\infty,r]}}
\newcommand{\chirr}{\one_{(-r,\,r)}}

\newcommand{\So}{S^{\circ}}

\newcommand{\Sr}{S_r}
\newcommand{\Srp}{\Sr^+}
\newcommand{\Srm}{\Sr^-}
\newcommand{\Sro}{\Sr^{\circ}}

\newcommand{\norm}[1]{\left\|#1\right\|}
\newcommand{\case}[1]{\begin{cases}#1\end{cases}}
\newcommand{\inv}{^{-1}}
\newcommand{\restr}[2]{\left. #1 \right|_{#2}}
\newcommand{\eps}{\varepsilon}

\newcommand{\matr}[1]{
\begin{pmatrix}
#1
\end{pmatrix}}

\newcommand{\smatr}[1]{\br{ \begin{smallmatrix}#1\end{smallmatrix} } }

\newcommand{\R}{\mathbb R}
\newcommand{\N}{\mathbb N}
\newcommand{\Z}{\mathbb Z}
\newcommand{\CC}{\mathbb C}

\DeclareMathOperator{\im}{Im}
\DeclareMathOperator{\Ker}{Ker}
\DeclareMathOperator{\Coker}{Coker}
\DeclareMathOperator{\Gr}{Gr}

\DeclareMathOperator{\Hom}{Hom}
\DeclareMathOperator{\End}{End}

\DeclareMathOperator{\Iso}{Iso}

\DeclareMathOperator{\Id}{Id}

\DeclareMathOperator{\supp}{supp}

\DeclareMathOperator{\dom}{dom}
\DeclareMathOperator{\sign}{sign}

\DeclareMathOperator{\W}{\mathcal{W}}

\DeclareMathOperator{\Proj}{\mathcal P}
\newcommand{\U}{\mathcal U}
\newcommand{\UK}{\U_K}

\newcommand{\B}{\mathcal B}
\newcommand{\Beu}{\B\eu}

\newcommand{\Reg}{\mathcal R}
\newcommand{\Rsa}{\Reg\sa}

\newcommand{\CRR}{\Reg_K}

\newcommand{\K}{\mathcal K}

\newcommand{\ZZ}{\mathcal{Z}}
\renewcommand{\S}{\mathcal{S}}
\newcommand{\Gs}{\mathcal{G}}

\renewcommand{\H}{\mathcal H}

\newcommand{\T}{\mathcal T}

\newcommand{\A}{\mathcal A}
\newcommand{\C}{\mathcal C}

\newcommand{\D}{\mathcal D}
\newcommand{\V}{\mathcal V}
\newcommand{\Y}{\mathcal Y}

\DeclareMathOperator{\ind}{ind}

\DeclareMathOperator{\inj}{inj}
\DeclareMathOperator{\dist}{dist}

\title{Spectral Sections} 
\author{Marina\,Prokhorova\thanks{Department of Mathematics, Technion -- Israel Institute of Technology}
\thanks{This work was partially supported by ISF grant no. 431/20}}
\date{}

\begin{document}

\maketitle

\uppskip
\begin{abstract}
The paper is devoted to the notion of a spectral section introduced by Melrose and Piazza.
In the first part of the paper we generalize results of Melrose and Piazza 
to arbitrary base spaces, not necessarily compact.
The second part contains a number of special cases, 
including cobordism theorems for families of Dirac type operators 
parametrized by a non-compact base space.
In the third part of the paper we investigate whether Riesz continuity is necessary 
for existence of a spectral section or a generalized spectral section.
In particular, we show that if a graph continuous family 
of regular self-adjoint operators with compact resolvents has a spectral section,
then the family is Riesz continuous.
\end{abstract}

\let\oldnumberline=\numberline
\renewcommand{\numberline}{\vspace{-19pt}\oldnumberline}
\addtocontents{toc}{\protect\renewcommand{\bfseries}{}}
\addtocontents{toc}{\vspace{-18pt}}

\uppskip

\tableofcontents

\section{Introduction}

In this paper we deal with families of regular (that is, closed and densely defined) 
linear operators with compact resolvents acting between Hilbert spaces
and parametri\-zed by points of an arbitrary topological space $X$.

Throughout the paper, \q{Hilbert space} always means a separable complex Hilbert space of infinite dimension;
\q{projection} always means an orthogonal projection, that is, a self-adjoint idempotent;
$\B(H,H')$ denotes the space of bounded linear operators from $H$ to $H'$ with the norm topology,
$\K(H,H')$ denotes the subspace of $\B(H,H')$ consisting of compact operators,
$\Bsa$ denotes the subspace of $\B(H)=\B(H,H)$ consisting of self-adjoint operators,
and $\Proj(H)$ denotes the subspace of $\Bsa$ consisting of projections. 
\index{1@Spaces!1B@$\B(H)$, bounded operators}
\index{1@Spaces!1B2sa@$\Bsa$, self-adjoint operators}
\index{1@Spaces!1K@$\K(H)$, compact operators}
\index{1@Spaces!1U@$\U(H)$, unitary operators}
\index{1@Spaces!1P@$\Proj(H)$, projections}

\sub{Spectral sections.}
\index{spectral section}
The notion of a spectral section was introduced by Melrose and Piazza in \cite{MP1},
in order to give a family version of a global boundary condition of Atiyah-Patody-Singer type.
It is convenient to split \cite[Definition 1]{MP1} into two parts, as we do below.

Let $A$ be a regular self-adjoint operator with compact resolvents. 
A projection $P$ is called an \textit{$r$-spectral section} for $A$ if
\begin{equation}\label{eq:ss-def00}
 Au = \lambda u \quad \Longrightarrow \quad
\case{ P u = u \;\text{ if } \lambda\geq r \\ 
       P u = 0 \;\text{ if } \lambda\leq -r }
\end{equation}
In other words,
\begin{equation}\label{eq:ss-def0}
 \one_{[r,+\infty)}(A) \leq P \leq \one_{(-r,+\infty)}(A).
\end{equation}
Here $\one_S$ denotes the characteristic function of the subset $S\subset\R$, 
and we use the natural partial order on the set of projections.
\index{2@Maps!1o@$\one$, characteristic function}
We also say that $P$ is a spectral section for $A$ with a \textit{cut-off parameter} $r$.

Let now $\A = (\A_x)_{\xX}$ be a family of regular self-adjoint operators with compact resolvents
and $r\colon X\to\R_+$ be a continuous function. 
A norm continuous family $P = (P_x)_{\xX}$ of projections 
is called an \textit{$r$-spectral section} for $\A$
if, for every $\xX$, the projection $P_x$ is an $r_x$-spectral section for the operator $\A_x$.
We also say that $P$ is a spectral section for $\A$ with a \textit{cut-off function} $r$.
\index{cut-off function}

Point out that the requirement for the family $(P_x)$ of projections 
and for the cut-off function $r$ to be continuous
is an essential part of the definition of a spectral section.

In \eqref{eq:ss-def00} we replaced strict inequalities 
used by Melrose and Piazza with non-strict ones,
since it simplifies statements of our results.
It influences only the cut-off function and does not change the notion of a spectral section.

\sub{Generalized spectral sections.}
\index{spectral section!generalized}
The notion of a generalized spectral section was introduced by Dai and Zhang in \cite{DZ1}
in order to cover both usual spectral sections and the Calder\'on projection. 
A projection $P$ is called a generalized spectral section for a self-adjoint operator $A$ 
if $P-\chip(A)$ is a compact operator.
A norm continuous family $P = (P_x)_{\xX}$ of projections 
is called a \textit{generalized spectral section} for a family $\A=(\A_x)$ of self-adjoint operators 
if $P_x$ is a generalized spectral section for $\A_x$ for every $\xX$. 
(We omit a requirement from \cite[Definition 2.1]{DZ1} 
for projections $P_x$ to be pseudodifferential, 
since we work in the general functional-analytic framework in the most part of the paper.)

\sub{Compact base spaces.}
Let $\A$ be a family of first order elliptic self-adjoint differential operators over a closed manifold
parametrized by points of a compact space $X$.
As was shown by Melrose and Piazza in \cite[Proposition 1]{MP1},
for such a family $\A$ the existence of a spectral section is equivalent 
to the vanishing of the index of $\A$ in $K^1(X)$. 

Recall that the \emph{Riesz topology} 
on the set $\Reg(H,H')$ of regular operators from $H$ to $H'$
is induced by the \emph{bounded transform} map 
\[ f\colon\Reg(H,H')\to\B(H,H'), \quad f(A) = A(1+A^*A)^{-1/2}, \] 
from the norm topology on $\B(H,H')$.
\index{topology!Riesz} 
If a family of elliptic differential operators over a closed manifold has continuously changing coefficients,
then the corresponding family of regular operators 
acting between the Hilbert spaces of square-integrable sections of corresponding vector bundles
is Riesz continuous.

The proof of Melrose and Piazza does not use the specifics of differential operators, 
so their result can be stated in a more general form:

\begin{prop}\label{prop:MP1}
Let $X$ be a compact space and
$\A = (\A_x)_{\xX}$ be a Riesz continuous family of self-adjoint regular operators with compact resolvents
acting on a Hilbert space $H$. 
Then the following two conditions are equivalent:
 \begin{enumerate}[topsep=-3pt, itemsep=0pt, parsep=3pt, partopsep=0pt]
	\item $\A$ has a spectral section.
	\item $\ind(\A) = 0 \in K^1(X)$.
\end{enumerate}
\end{prop}

\upskip
Melrose and Piazza also proved a $\Z_2$-graded analog of this result in \cite[Proposition 2]{MP2}.
The proofs of these results in \cite{MP1,MP2} depend crucially on the base space being compact.

\sub{Arbitrary base spaces.}
The aim of the first part of this paper is to generalize aforementioned results of Melrose and Piazza 
to arbitrary base spaces, not necessarily compact.
In particular, we prove the following generalization of Proposition \ref{prop:MP1}.

\textbf{Theorem \ref{thm:ss}.}
\textit{ Let $X$ be an arbitrary topological space
  and $\A = (\A_x)_{\xX}$ be a Riesz continuous family of 
	regular self-adjoint operators with compact resolvents acting on a Hilbert space $H$.
  Then the following three conditions are equivalent: 
  \begin{enumerate}[topsep=-3pt, itemsep=0pt, parsep=3pt, partopsep=0pt]
	   \item $\A$ has a spectral section.
	   \item $\A$ has a generalized spectral section.
	   \item $\A$ is homotopic to a family of invertible operators.
  \end{enumerate}
  If $P$ is a generalized spectral section for $\A$, then for every $\eps>0$
  a spectral section $Q$ for $\A$ can be chosen so that
  $\norm{Q-P}_{\infty} < \eps$ and $Q$ is homotopic to $P$ as a generalized spectral section.
}

\medskip
We also prove a $\Z_2$-graded version of this result, see Theorem \ref{thm:ss1}. 
It deals with $\Cl(1)$ spectral sections for odd self-adjoint operators
and generalizes \cite[Pro\-po\-sition 2]{MP2} in the same manner as 
Theorem \ref{thm:ss} generalizes \cite[Proposition 1]{MP1}.

\sub{Trivializing operators of finite rank.}
Melrose and Piazza showed in \cite[Lemma 8]{MP1} 
that if a self-adjoint family $\A$ over a compact base space has a spectral section $P$, 
then there is a \q{trivializing} family $\C$ of self-adjoint finite rank operators
such that the \q{trivialized} family $\A'=\A+\C$ consists of invertible operators and
$P$ is the family of positive spectral projections for $\A'$.
They also proved a $\Z_2$-graded analog in \cite[Lemma 1]{MP2}.
We generalize both these results to arbitrary base spaces in Sections \ref{sec:tr-op}--\ref{sec:Cl1}.

\sub{Applications.}
In Part \ref{part:appl} we present a variety of special cases 
illustrating how the results of Part \ref{part:Riesz} can be used.
In particular, in Section \ref{sec:cob} 
we generalize the famous Cobordism Theorem for Dirac operators 
to arbitrary, not necessarily compact, base spaces. 
We show in Theorem \ref{thm:cobD} that for an arbitrary family of Dirac type operators 
on a compact manifold with boundary,
the family of symmetrized boundary operators has a spectral section.
$\Z/2$-graded case of this result is given by Theorem \ref{thm:cobD1}.
In our cobordism theorems, we do not require the product form of the operator near the boundary,
neither we require the product form of the metric.

\sub{Graph continuous families.}
The Riesz topology on the set of regular operators 
is well suited for the theory of differential operators on closed manifolds.
However, it is not quite adequate for differential operators on manifolds with boundary:
it is unknown, except for several special cases, whether families of regular operators 
defined by boundary value problems are Riesz continuous.
In the context of boundary value problems, the better suited topology is the \textit{graph topology}.

The graph topology on $\Reg(H,H')$ is induced by the inclusion of $\Reg(H,H')$ to $\Proj(H\oplus H')$ 
taking a regular operator to the orthogonal projection onto its graph.
\index{topology!graph}
A family of elliptic operators and boundary conditions 
with continuously varying coefficients leads to a graph continuous family of regular operators
\cite[Appendix A.5]{Pr17}.

A spectral section of a Riesz continuous family always exists locally.
As we saw above, there is a topological obstruction for existence of a global spectral section, 
which in the case of a compact base space takes value in the $K^1$-group of the base.
A \textit{graph continuous} family, however, may admit \textit{no spectral section even locally}.
In fact, Riesz continuity is \textit{necessary} for a local existence of a spectral section,
as the following result shows.

\textbf{Theorem \ref{thm:ss2Riesz}.}
\textit{Let $\A = (\A_x)_{\xX}$ be a graph continuous family of 
	regular self-adjoint operators with compact resolvents acting on a Hilbert space $H$.
  Suppose that $\A$ admits a spectral section.
  Then $\A$ is Riesz continuous.
}

The situation with \textit{generalized} spectral sections for graph continuous families is more complicated,
as we show in Section \ref{sec:gt-gss}.
On the one hand, a generalized spectral section does not necessarily exist
even when a base space is an interval and/or a family consists of invertible operators.
On the other hand, existence of a generalized spectral section does not imply Riesz continuity.

\sub{Acknowledgments.}
I am grateful to N.V.~Ivanov 
for useful remarks and suggestions.

\section{Preliminaries}

In this section we recall some basic facts about regular operators
(for more detailed exposition see, for example, \cite{L04, BBLP, Kato})
and introduce some designations that are used in the rest of the paper.

\sub{Regular operators.}
\index{operators!regular}
\index{1@Spaces!1R0@$\Reg(H)$, regular operators}
An unbounded operator $A$ from $H$ to $H'$ is a linear operator 
defined on a subspace $\D$ of $H$ and taking values in $H'$;
the subspace $\D$ is called the domain of $A$ and is denoted by $\dom(A)$.
An unbounded operator $A$ is called closed if its graph is closed in $H\oplus H'$ 
and densely defined if its domain is dense in $H$.
It is called \textit{regular} if it is closed and densely defined.
Let $\Reg(H,H')$ denote the set of all regular operators from $H$ to $H'$,
and let $\Reg(H) = \Reg(H,H)$.

The \emph{adjoint operator} of an operator $A\in\Reg(H,H')$ 
is an unbounded operator $A^*$ from $H'$ to $H$ with the domain 
$$\dom(A^*) = \sett{u\in H'}{\text{ there exists } v\in H \text{ such that } 
   \bra{Aw,u}=\bra{w,v} \text{ for all } w\in H}.$$
For $u\in\dom(A^*)$ such an element $v$ is unique and $A^*u=v$ by definition.
The adjoint of a regular operator is itself a regular operator.

An operator $A\in\Reg(H)$ is called \emph{self-adjoint} if $A^*=A$ (in particular, $\dom(A^*) = \dom(A)$).
Let $\Reg\sa(H)\subset \Reg(H)$ be the subset of self-adjoint regular operators.
\index{operators!regular!self-adjoint}
\index{1@Spaces!1R2sa@$\Reg\sa(H)$, self-adjoint regular operators}

\sub{Operators with compact resolvents.}
For a regular operator $A\in\Reg(H,H')$, 
the operator $1+A^*A$ is regular, self-adjoint, and has dense range. 
Its densely defined inverse $(1+A^*A)\inv$ is bounded 
and hence can be extended to a bounded operator defined on the whole $H$.
A regular operator $A\in\Reg(H,H')$ is said to \textit{have compact resolvents} 
if $(1+A^*A)\inv\in\K(H)$ and $(1+AA^*)\inv\in\K(H')$.
We denote by $\CRR(H)$ the subset of $\Reg(H)$ consisting of regular operators with compact resolvents.
\index{operators!regular!with compact resolvents}
\index{1@Spaces!1R2K@$\CRR(H)$, regular operators with compact resolvents}

Let $\CRR\sa(H) = \Reg\sa(H)\cap\CRR(H)$ be the subset of $\Reg(H)$ 
consisting of self-adjoint operators with compact resolvents.
Equivalently, a self-adjoint regular operator $A$ is an operator with compact resolvents
if $(A+i)\inv$ is a compact operator.
Such an operator has a discrete real spectrum.

\sub{Bounded transform.}
\index{bounded transform}
\index{2@Maps!2f@$f$, bounded transform}
The \emph{bounded transform} (or the Riesz map) 
$A\mapsto f(A) = A(1+A^*A)^{-1/2}$
defines the inclusion of the set $\Reg(H,H')$ of regular operators 
to the unit ball in the space $\B(H,H')$ of bounded operators, with the image
\begin{equation*}
	f(\Reg(H,H')) = \sett{ a\in\B(H,H')}{\norm{a}\leq 1 \text{ and } 1-a^*a \text{ is injective} }.
\end{equation*}
The inverse map is given by the formula $f\inv(a) = a(1-a^*a)^{-1/2}$.

If $A$ is self-adjoint, then so is $f(A)$.
If $A$ has compact resolvents, then $a=f(A)$ is essentially unitary 
(that is, both $1-a^*a$ and $1-aa^*$ are compact operators).

\sub{Essentially unitary operators.}
\index{operators!bounded!essentially unitary}
\index{1@Spaces!1B3eu@$\B\eu(H)$, essentially unitary operators}
Let $\B\eu(H,H')$ be the subspace of $\B(H,H')$ consisting of essentially unitary operators.
The bounded transform takes $\CRR(H,H')$ to
\begin{equation*}
	f\br{\CRR(H,H')} = \sett{ a\in\B\eu(H,H')}{\norm{a}\leq 1 \text{ and } 1-a^*a \text{ is injective} }.
\end{equation*}
Let $\Beusa$ be the the subspace of $\B(H,H')$ consisting of self-adjoint essentially unitary operators.

In Section \ref{sec:gss} and the first part of Section \ref{sec:Cl1} we prove our results simultaneously 
both for regular operators with compact resolvents and for essentially unitary operators.

The notion of a generalized spectral section given in the Introduction 
works as well for self-adjoint essentially unitary operators.
\index{generalized spectral section}
Equivalently, a projection $P$ is a generalized spectral section for $a\in\Beusa$
if $(2P-1)-a\in\K(H)$
(since $\chip(a)$ is a compact deformation of $(a+1)/2$ for $a\in\Beusa$).

We call a projection $P$ an $r$-spectral section for $a\in\Beusa$ and $r\in(0,1)$
\index{spectral section}
if $\one_{[r,+\infty)}(a) \leq P \leq \one_{(-r,+\infty)}(a)$.
Let $r\colon X\to(0,1)$ be a continuous function; 
we call a norm continuous family of projections $P=(P_x)_{\xX}$ 
an $r$-spectral section (or simply a spectral section)
for a family $a=(a_x)$ of self-adjoint essentially unitary operators
if $P_x$ is an $r_x$-spectral section for $a_x$ for every $\xX$.

A family $P=(P_x)$ of projections is a generalized spectral section, resp. $r$-spectral section
for a family $\A=(\A_x)$ of self-adjoint regular operators with compact resolvents
if and only if $P$ is a generalized spectral section, resp. $(f\circ r)$-spectral section
for the family $f\circ\A$ of self-adjoint essentially unitary operators.

\sub{Two topologies on $\Reg(H,H')$.}
There are several natural topologies on the set of regular operators.
The two most useful of them are the Riesz topology and the graph topology.
They are induced by the inclusions of $\Reg(H,H')$ to $\B(H,H')$ and to $\Proj(H\oplus H')$, respectively;
see Introduction for details.

We will always specify what topology (Riesz or graph) on $\Reg(H,H')$ we consider.
We will write $^r\Reg(H,H')$ or $^g\Reg(H,H')$
for the space of regular operators with Riesz or graph topology, respectively.
\index{1@Spaces!1R1@$^r\Reg(H)$, $^g\Reg(H)$}
Alternatively, we will write \q{Riesz continuous} or \q{graph continuous} 
for maps from or to $\Reg(H,H')$ and for families of regular operators.

On the subset $\B(H,H')\subset \Reg(H,H')$ both Riesz and graph topology coincide with the usual norm topology. 
Therefore, we always consider $\B(H,H')$ as equipped with the norm topology.

\sub{Spectral projections.}
The spectrum of an operator $A\in\CRR\sa(H)$ is discrete and real. 
For real numbers $a<b$, the spectral projection $\one_{[a,b]}(A)$ is defined as 
$\one_{[f(a),\,f(b)]}(f(A))$;
its range is the subspace of $H$ spanned by eigenvectors of $A$ with eigenvalues in the interval $[a,b]$.
Similarly, $\one_{[a,+\infty)}(A)$ is defined as $\one_{[f(a),1]}(f(A))$ and
$\one_{(-\infty,a]}(A)$ is defined as $\one_{[-1,f(a)]}(f(A))$.
The spectral projections for semi-open and open intervals are defined in the same manner.

Let  $\Res(A)$ denote the resolvent set of $A$.
For a compact subspace $K\subset\R$,
the subset $$V_K = \sett{A\in\CRR\sa(H)}{K\subset\Res(A)}$$
is open in both Riesz and graph topology on $\CRR\sa(H)$.

The map $V_{\set{a,b}}\to\Proj(H)$ taking $A$ to $\one_{[a,b]}(A)$
is both Riesz-to-norm and graph-to-norm continuous. 
However, the spectral projection maps $V_{a}\to\Proj(H)$ corresponding to unbounded intervals,
$A\mapsto\one_{(-\infty,a]}(A)$ and $A\mapsto\one_{[a,+\infty)}(A)$,
are only Riesz-to-norm continuous, but not graph-to-norm continuous.
This is the major difference between the two topologies in our context.

\addtocontents{toc}{\vspace{-5pt}}

\part{Riesz continuous families}\label{part:Riesz}

\addtocontents{toc}{\vspace{-12pt}}

Throughout this part, all families of regular operators are supposed to be Riesz continuous.

\upskip
\section{Generalized spectral sections}\label{sec:gss}

\upskip
\sub{Homotopy Lifting Property.}
Recall that a continuous map $\Y\to\ZZ$ between topological spaces
is said to have the Homotopy Lifting Property for a space $X$ 
if for every homotopy $h\colon X\times[0,1]\to\ZZ$, 
every lifting $\tilde{h}_0\colon X\times\set{0}\to\Y$ of $h_0$
can be continued to a lifting $\tilde{h}\colon X\times[0,1]\to\Y$ of $h$.

In this section the base space $\ZZ$ is either $^r\CRRsa$ or $\Beusa$.
\index{1@Spaces!2Z@$\ZZ$}
For the most part of the section, our reasoning  works for both these cases simultaneously.

Let	$\ZZ$ be one of these two spaces.
Let $I$ denote the range of cut-off parameters for $\ZZ$,
that is, $I=\R_+$ if $\ZZ=\CRRsa$ and $I=(0,1)$ if $\ZZ=\Beusa$.
\index{5I@$I$, range of cut-off parameters}

Let $\Gs$ denote the subspace of $\ZZ\times\Proj(H)$ consisting of pairs 
$(A,P)$ such that $P$ is a generalized spectral section for $A$.
We consider $\Gs$ as the total space of 
a fiber bundle over $\ZZ$, with the projection $\pig\colon\Gs\to \ZZ$ taking $(A,P)$ to $A$.
\index{3@Bundles and bundle maps!2piG@$\pig\colon\Gs\to\ZZ$}

\begin{thm}\label{thm:Omega}
	The fiber bundle $\pig\colon\Gs\to \ZZ$ is locally trivial 
	and has the Homotopy Lifting Property for all spaces.
\end{thm}

\proof 
The family $(\ZZ_r)_{r\in I}$, $\ZZ_r = \sett{A\in \ZZ }{ r\in\Res(A)}$, 
is an open covering of $\ZZ$.
For every $r\in I$, the map $\Srp\colon \ZZ_r\to\Proj(H)$ 
defined by the formula $\Srp(A) = \chirp(A)$ is continuous.
The restriction of $\Gs$ to $\ZZ_r$ is  
\begin{equation}\label{eq:OZl}
  \restr{\Gs}{\ZZ_r} = \sett{(A,P)}{A\in \ZZ_r \text{ and } P-\Srp(A)\in\K(H) }. 
\end{equation}
Our next goal is to trivialize $\Srp$ locally.
To this end, take an open covering $(V_Q)_{Q\in\Proj(H)}$ of $\Proj(H)$,
where $V_Q = \sett{P\in\Proj(H)}{\norm{P-Q}<1 }$.
By \cite[Proposition 5.2.6]{WO} for every $Q\in\Proj(H)$ 
there is a continuous map $g_Q\colon V_Q\to\U(H)$
such that 
\begin{equation}\label{eq:GQ}
	P = g_Q(P)\, Q\, (g_Q(P))^* \; \text{ for every } P\in V_Q.
\end{equation}
The inverse images $\ZZ_{r,\,Q} = (\Srp)\inv(V_Q)\subset \ZZ_r$, 
with $r$ running $I$ and $Q$ running $\Proj(H)$, 
form an open covering $(\ZZ_{r,\,Q})$ of $\ZZ$.
We claim that the restriction $\Gs_{r,\,Q}$ of $\Gs$ to $\ZZ_{r,\,Q}$ 
is a trivial bundle with the fiber 
\[ F_Q = \sett{P\in\Proj(H)}{P-Q\in\K(H) }. \] 
Indeed, fix an arbitrary pair $(r,Q)$ and consider the map 
$g = g_Q\circ \Srp\colon \ZZ_{r,\,Q}\to\U(H)$.
It follows from \eqref{eq:GQ} that 
$(g(A))^*\, \Srp(A)\, g(A) = Q$
does not depend on $A\in\ZZ_{r,\,Q}$.
Together with \eqref{eq:OZl} this implies that the map 
\[ \Phi\colon\Gs_{r,\,Q}\to\ZZ_{r,\,Q}\times F_Q, \quad \Phi(A,P) = (A,\, g(A)^*\, P \, g(A)) \]
is a trivializing bundle isomorphism, which proves the claim.

For both $\ZZ = {^r\CRRsa}$ and $\ZZ=\Beusa$
the topology of $\ZZ$ is induced by the embedding of $\ZZ$ to $\B(H)$,
so $\ZZ$ is a metric space. 
This implies paracompactness of $\ZZ$ \cite[Corollary 1]{St}.
Thus $\pig$ is a locally trivial fiber bundle with a paracompact base space.
By \cite[Uniformization Theorem]{Hur}, 
$\pig$ has the Homotopy Lifting Property for all spaces.
This completes the proof of the proposition.
\endproof

\sub{Generalized spectral sections.}
Using Theorem \ref{thm:Omega}, we now can prove the following result.

\begin{thm}\label{thm:gss-inv}
  Let $X$ be an arbitrary topological space and
	$\ZZ$ be either $^r\CRRsa$ or $\Beusa$.
	Let $\A\colon X\to \ZZ$ be a continuous map.
	Then the following two conditions are equivalent:
\begin{enumerate}[topsep=-3pt, itemsep=0pt, parsep=3pt, partopsep=0pt]
	\item $\A$ has a generalized spectral section.
	\item $\A$ is homotopic (as a map from $X$ to $\ZZ$) to a family of invertible operators.
\end{enumerate}
\end{thm}

\proof
($2\Rightarrow 1$)
Let $h\colon X\times[0,1]\to\ZZ$ be a homotopy between $\A=h_1$ and an invertible family $h_0$.
Since $h_0$ is invertible, it has a spectral section $P_0(x) = \chip(h_0(x))$. 
Then the map $\tilde{h}_0\colon X\to\Gs$ given by the formula $\tilde{h}_0(x) = (h_0(x),P_0(x))$ covers $h_0$.
By Theorem \ref{thm:Omega}, 
$\tilde{h}_0$ can be continued to a map $\tilde{h}\colon X\times[0,1]\to\Gs$ covering $h$.
Restriction of $\tilde{h}$ to $X\times\set{1}$ gives the map $\tilde{h}_1\colon X\to\Gs$ covering $\A$.
The composition of $\tilde{h}_1$ with the natural projection $\Gs\to\Proj(H)$ 
is a generalized spectral section for $A$.

($1\Rightarrow 2$)
Let $P\colon X\to \Proj(H)$ be a generalized spectral section for $\A\colon X\to \ZZ$.
Then $T\colon X\to\B(H)$, $T(x) = 2P(x)-1$, is a continuous family of symmetries (that is, self-adjoint unitaries).

1. Consider first the case $\ZZ=\Beusa$.
Then $\A(x)-T(x)\in\K(H)$ for every $\xX$.
Therefore, 
\begin{equation*}\label{hBeu}
	h\colon X\times[0,1]\to\Beusa, \quad h_t(x) = (1-t)\A(x)+tT(x)
\end{equation*}
is a homotopy from $h_0=\A$ to $h_1=T$, with $T(x)$ invertible for every $\xX$.

2. Let now $\ZZ = {^r\CRRsa}$.
The composition $a = f\circ\A\colon X\to\Beusa$ 
is continuous and $a(x)-T(x)\in\K(H)$ for every $\xX$.
One is tempted to apply $f\inv$ to the linear homotopy between $a$ and $T$, as above.
But the image of $f$, 
\begin{equation}\label{eq:fCRR}
	f\br{\CRRsa} = \sett{b\in\Beusa}{\norm{b}\leq 1 \text{ and } 1-b^2 \text{ is injective }},
\end{equation}
does not contain $T(x)$, so this naive idea does not work.
To fix it, we replace $T$ by its compact deformation $T'$ lying in the image of $f$.
Let us fix a strictly positive compact operator $K\in\K(H)$ of norm less than $1$.
For example, one can identify $H$ with $l^2(\N)$ and take the diagonal operator 
$K=\diag(\frac{1}{2},\frac{1}{3},\frac{1}{4},\ldots)$.
Then $T'(x) = (1-K)T(x)(1-K)$ is self-adjoint and invertible, $\norm{T'(x)}\leq 1$, and $T'(x)-T(x)\in\K(H)$.
Let $a_t(x) = (1-t)a(x)+tT'(x)$ be the linear homotopy between $a_0=a$ and $a_1=T'$.
By definition, $a_0(x) = f(\A(x))$ lies in the image of $f$.
For every $t\in(0,1]$, $\xX$, and $\xi\in H\setminus\set{0}$ we get $\norm{a_t(x)\xi} < \norm{\xi}$, 
so $1-a_t(x)^2$ is injective and thus $a_t(x)$ lies in the image of $f$.
Applying $f\inv$ to $a_t(x)$, we obtain a homotopy $h\colon X\times[0,1]\to\CRRsa$, 
$h_t(x) = f\inv(a_t(x))$, connecting $h_0=\A$ with an invertible family $f\inv(T')$.
This completes the proof of the theorem.
\endproof

\section{Spectral sections}\label{sec:ss}

In this section $\ZZ$ always denotes the space $\CRRsa$ equipped with the Riesz topology.

\sub{Fiber homotopy equivalence.}
Let $\S$ be the subspace of $\ZZ\times\Proj(H)\times\R_+$ consisting of triples 
$(A,P,r)$ such that $P$ is an $r$-spectral section for $A$.
We consider $\S$ as the total space of 
a fiber bundle over $\ZZ$, with the projection $\pis\colon\S\to\ZZ$ taking $(A,P,r)$ to $A$.
\index{3@Bundles and bundle maps!2piS@$\pis\colon\S\to\ZZ$}

\begin{thm}\label{thm:ss-gs}
  The bundle map $\iota\colon\S\to\Gs$ taking $(A,P,r)$ to $(A,P)$
	is a fiber homotopy equivalence.
	Moreover, for every $\eps>0$, a fiber homotopy inverse bundle map $\varphi = \varphi_{\eps} \colon\Gs\to\S$
	can be chosen so that $\norm{Q-P} < \eps$ for every $(A,P)\in\Gs$, 
	$(A,Q) = \iota\circ\varphi(A,P)$.\index{3@Bundles and bundle maps!3i@$\iota\colon\S\to\Gs$}
\index{3@Bundles and bundle maps!4phi@$\varphi_{\eps}\colon\Gs\to\S$}
\end{thm}

As an immediate corollary of this theorem we get the following result.

\begin{cor}
  The fiber bundle $\pis\colon\S\to\ZZ$ has the Weak Homotopy Lifting Property for all spaces
	(that is, for any homotopy $h\colon X\times[0,1]\to\ZZ$ 
and for any lifting $\tilde{h}_0\colon X\times\set{0}\to\S$ of $h_0$,
there is a lifting $X\times[0,1]\to\S$ of $h$ whose restriction to $X\times\set{0}$ 
is vertically homotopic to $\tilde{h}_0$).
\end{cor}

\prooff{ of Theorem \ref{thm:ss-gs}}
Without restriction of generality, we can suppose that $\eps<1$.
Fix $\delta\in(0,\frac{\eps}{2})$.

\sloppy{
Let $(A,P)$ be an arbitrary element of $\Gs$.
For $r\geq 0$, define the spectral projections 
\begin{equation}\label{eq:Sr}
		\Srm(A) = \chirm(A), \quad
		\Sro(A) = \chirr(A), \; \text{ and } \;
		\Srp(A) = \chirp(A).
\end{equation}
We approximate $P$ by bounded self-adjoint operators	
\begin{equation}\label{eq:Tr}
	T_r(A,P) = \Srp(A) + \Sro(A)\,P\,\Sro(A).
\end{equation}
The map $T_r\colon\Gs\to\B(H)$ so defined is continuous on the open subset 
\[ \sett{(A,P)\in\Gs}{\pm r\in\Res(A)} \subset\Gs. \]
Equality \eqref{eq:Tr} can be written equivalently as 
\begin{equation}\label{eq:Tr0}
	T_r(A,P) = S^+_0(A) + \Sro(A)\br{P-S^+_0(A)}\Sro(A).
\end{equation}
Consider the family 
\begin{equation*}\label{eq:Ur}
	\U = (U_r)_{r>0}, \quad
	U_r = \sett{(A,P)\in\Gs}{\pm r\in\Res(A) \text{ and } \norm{T_r(A,P)-P}<\delta }
\end{equation*}
of open subsets of $\Gs$.
We claim that $\U$ is an open covering of $\Gs$.
Indeed, let $(A,P)$ be an arbitrary point of $\Gs$.
Since $P$ is a generalized spectral section for $A$,
the difference $P-S^+_0(A)$ is a compact operator.
The net $\set{\Sro(A)}_{r>0}$ is an approximate unit for $\K(H)$. 
Therefore, the second summand in the right hand side of \eqref{eq:Tr0} has a limit $P-S^+_0(A)$,
and $T_r(A,P)-P\to 0$ as $r\to +\infty$.
It follows that $(A,P)$ lies in $U_r$ for $r$ large enough,
and thus $\U$ covers $\Gs$.
}

The topology of $\ZZ\times\Proj(H)$ is induced by its embedding to $\B(H)\times\B(H)$,
so $\ZZ\times\Proj(H)$, as well as its subspace $\Gs$, is a  metric space. 
It follows from \cite[Corollary 1]{St} that $\Gs$ is paracompact.
Hence there is a partition of unity $(u_i)$ subordinated to $\U$,
with $\supp(u_i)\subset U_{r_i}$.
(By a partition of unity we always mean \textit{locally finite} partition.)
For every $(A,P)\in\Gs$ we define a bounded self-adjoint operator $T(A,P)$ by the formula
\begin{equation}\label{eq:T}
	T(A,P) = \sum u_i(A,P)\cdot T_{r_i}(A,P).
\end{equation}
The restriction of $T_r$ to $U_r$ is continuous, 
so all the summands in \eqref{eq:T} are continuous,
and thus $T$ itself is continuous as a map from $\Gs$ to $\B(H)$.

If $u_i(A,P)\neq 0$, then $(A,P)\in U_{r_i}$ and $\norm{T_{r_i}(A,P)-P} < \delta$.
Therefore, 
\begin{equation}\label{eq:TP}
	\norm{T(A,P)-P} < \delta \text{ for every } (A,P)\in\Gs.
\end{equation}
The spectrum of $P$ is $\set{0,1}$, 
so the last inequality implies that the spectrum of $T(A,P)$ lies in 
$\Lambda = [-\delta,\,\delta]\cup[1-\delta,\,1+\delta]$.
Since $\delta<\frac{1}{2}$, these two intervals are disjoint. 
The function $\one_{[\frac{1}{2},\,+\infty)}$ is continuous (and even smooth) on $\Lambda$,
so the projection
\begin{equation*}\label{eq:Q}
	Q(A,P) = \one_{[\frac{1}{2},\,+\infty)}(T(A,P))
\end{equation*}
depends continuously on $(A,P)$.
Moreover, $\norm{Q(A,P)-T(A,P)} \leq \delta$, 
which together with \eqref{eq:TP} provides an estimate
$\norm{Q(A,P)-P} < \eps$ for every $(A,P)\in\Gs$.

Let $(A,P)\in\Gs$ and $r=\max\sett{r_i}{u_i(A,P)\neq 0}$. 
Then 
 $T(A,P) = 0\oplus T^{\circ}(A,P)\oplus 1$ and thus $Q(A,P) = 0\oplus Q^{\circ}(A,P)\oplus 1$
with respect to the orthogonal decomposition 
\begin{equation*}\label{eq:decH}
  H = H^-\oplus H^{\circ} \oplus H^+, \quad \text{ where } 
	H^{\alpha} = \im(S^{\alpha}_r(A)) \text{ for } \alpha\in\set{+,-,\circ}.
\end{equation*}
In other words, $Q(A,P)$ is an $r$-spectral section for $A$. 

For every point $(A,P)\in\Gs$ choose its neighbourhood $V_{A,P}\subset\Gs$ 
intersecting only a finite number of inverse images $u_i\inv(0,1]$.
Let $(v_j)$ be a partition of unity subordinated to the covering $\V = (V_{A,P})$ of $\Gs$.
Let $$R_j = \max\sett{r_i}{u_i\inv(0,1]\cap v_j\inv(0,1] \neq\emptyset}.$$
We define a continuous map $R\colon \Gs\to \R_+$ by the formula 
\begin{equation*}\label{eq:R}
	R(A,P) = \sum v_j(A,P)\cdot R_j.
\end{equation*}
Then $R(A,P) \geq \max\sett{r_i}{u_i(A,P)\neq 0}$ for every $(A,P)\in\Gs$,
so $Q(A,P)$ is an $R(A,P)$-spectral section for $A$.

Finally, we define the map $\varphi\colon\Gs\to\S$ by the formula $\varphi(A,P) = (A, Q(A,P), R(A,P))$.

It remains to show that $\varphi$ is fiber homotopy inverse to $\iota$.
This follows immediately from the following lemma.

\begin{lem}\label{lem:PQR}
Let $\varphi\colon\Gs\to\S$ be a continuous map, $\varphi(A,P) = (A, Q(A,P), R(A,P))$,
such that $\norm{Q(A,P)-P} < 1$ for every $(A,P)\in\Gs$.
Then $\varphi$ is fiber homotopy inverse to $\iota$.
\end{lem}

\proof
The formula 
\begin{equation*}\label{eq:PQt}
  [P_0,P_1]_t = \one_{[\frac{1}{2},\,+\infty)}\br{(1-t)P_0+tP_1}  
\end{equation*}
determines a continuous map 
\[ \Proj(H)^2\times[0,1] \;\supset\; \sett{(P_0,P_1,t)}{\norm{P_0-P_1}<1} \to \Proj(H), \]
since $(1-t)P_0+tP_1$ lies on the distance less than $1/2$ either from $P_0$ or from $P_1$
and thus its spectrum is contained in $[0,1/2)\cup(1/2,1]$.
Clearly, $[P_0,P_1]_t = P_t$ for $t=0,1$.
If both $P_0$ and $P_1$ are generalized spectral sections (respectively $r$-spectral sections) for $A$, 
then $[P_0,P_1]_t$ is also a generalized spectral section (respectively $r$-spectral section) for $A$.

We will write $Q$ and $R$ instead of $Q(A,P)$ and $R(A,P)$ for brevity. 
A bundle homotopy between $\Id_{\Gs}$ and $\iota\circ\varphi$ is given by the formula $h_t(A,P) = (A,[P,Q]_t)$.
To construct a homotopy between $\Id_{\S}$ and $\varphi\circ\iota$, 
we use three auxiliary homotopies.
The first homotopy $h'_t(A,P,r) = (A,P,r+tR)$ connects $(A,P,r)$ with $(A,P,r+R)$.
The second one $h''_t(A,P,r) = (A,[P,Q]_t,r+R)$ 
connects $(A,P,r+R)$ with $(A,Q,r+R)$.
The third one $h'''_t(A,P,r) = (A,Q,(1-t)r+R)$ connects $(A,Q,r+R)$ with $\varphi\circ\iota(A,P,r) = (A,Q,R)$.
The concatenation of $h'$, $h''$, and $h'''$ is a desired bundle homotopy 
between $\Id_{\S}$ and $\varphi\circ\iota$.
This completes the proof of the lemma and the theorem.
\endproof

\sub{Spectral sections.}
Using Theorem \ref{thm:ss-gs}, we immediately obtain the following result.

\begin{thm}\label{thm:ss}
  Let $X$ be an arbitrary topological space and $\A\colon X\to\CRRsa$ be a Riesz continuous map.
  Then the following three conditions are equivalent:
\begin{enumerate}[topsep=-3pt, itemsep=0pt, parsep=3pt, partopsep=0pt]
	\item $\A$ has a spectral section.
	\item $\A$ has a generalized spectral section.
	\item $\A$ is homotopic, via a Riesz continuous homotopy $X\times[0,1]\to\CRRsa$, 
	             to a family of invertible operators.
\end{enumerate}
		      If $P$ is a generalized spectral section for $\A$ and $\eps>0$, then
	        a spectral section $Q$ for $\A$ can be chosen so that
	        $\norm{Q-P}_{\infty} < \eps$ and $Q$ is homotopic to $P$ as a generalized spectral section.
\end{thm}

\proof
($1\Rightarrow 2$) is trivial.

($2\Leftrightarrow 3$) follows from Theorem \ref{thm:gss-inv}.

($2\Rightarrow 1$) and the last part of the theorem follow from Theorem \ref{thm:ss-gs}.
Indeed, if $\eps>0$ and $P$ is a generalized spectral section for $\A$,
then the map $x\mapsto\varphi_{\eps}(\A_x,P_x) =: (\A_x,Q_x,r_x)$ 
defines an $r$-spectral section $Q$ for $\A$.
Moreover, $\norm{Q_x-P_x} < \eps$ for every $\xX$.
Since $\iota\circ\varphi_{\eps}$ and $\Id_{\Gs}$ are vertically homotopic,
$Q$ and $P$ are homotopic as generalized spectral sections.
This completes the proof of the theorem.
\endproof

\begin{rem}\label{rem:large}  
For a self-adjoint family $\A_x$ parametrized by points of a compact space $X$, 
Melrose and Piazza showed in \cite[Proposition 2]{MP1}
that if the set of spectral sections for $\A$ is non-empty,
then it contains \q{arbitrary small} and \q{arbitrary large} spectral sections, 
in the following sense:
for every given $s\in\R$, there are spectral sections $P$ and $Q$ such that 
$P_x\leq S_x\leq Q_x$ for every $\xX$, where
	$S_x = \one_{[s,+\infty)}(\A_x)$.
This property is no longer true in the general case of a non-compact base space,
as the following simple example shows.
\end{rem}

\begin{exm}\label{exm:A+x}
Fix $A\in\CRR\sa(H)$ and consider the family $\A_x=A+x$ parametrized by real numbers $x\in X=\R$.
The family $\A$ admits a spectral section.
Indeed, the constant map taking every $\xX$ to $\chip(A)$ is an $r$-spectral section for $\A$,
where $r\colon X\to\R_+$ is an arbitrary function satisfying condition $r(x)>|x|$.

Suppose now that $A$ is unbounded from above. 
Then, for any given $s\in\R$, 
$\A$ has no spectral section $P$ dominated by the family $S=(S_x)$ as above.
Indeed, suppose that $P\colon X\to\Proj(H)$ is such a spectral section.
Then $P_x \leq \one_{[s,+\infty)}(A+x) = \one_{[s-x,\,+\infty)}(A)$.
For $x\leq s$, let $P'_x$ be the restriction of $P_x$ to the range $H'$ of $\one_{[0,+\infty)}(A)$.
Then $(P'_x)$ is a continuous one-parameter family of projections in $H'$
parametrized by points of the ray $(-\infty,s]$.
The kernels of $P'_x$ are finite-dimensional;
by continuity, their dimensions should be independent of $x$.
On the other hand, the dimension of $\Ker(P'_x)$ is bounded from below 
by the rank of $\one_{[0,\,s-x)}(A)$, which goes to infinity as $x\to -\infty$.
This contradiction shows that such a spectral section $P$ does not exist.

Similar argument shows that if $A$ is unbounded from below, 
then $\A$ has no spectral section dominating the family $S$.
\end{exm}

\section{Trivializing operators}\label{sec:tr-op}

Melrose and Piazza showed in \cite[Lemma 8]{MP1} 
that if a self-adjoint family $\A$ over a compact base space 
admits a spectral section $P$, 
then $\A$ admits a finite rank correction $\C$ to an invertible family $\A' = \A+\C$
such that $P$ is the family of positive spectral projections for $\A'$.
They also proved a $\Z_2$-graded analog of this result in \cite[Lemma 1]{MP2}.
Their proofs are based on the existence of 
\q{arbitrarily small} and \q{arbitrarily large} spectral sections in the sense of Remark \ref{rem:large}.
However, for a non-compact base space there may be no such spectral sections,
as Example \ref{exm:A+x} shows.
In this section we generalize \cite[Lemma 8]{MP1} to arbitrary base spaces using a different method.
A generalization of \cite[Lemma 1]{MP2} will be given in the next section.

\sub{Trivializing operators.}
Following \cite[Definition 2.11]{LP},
we use the term \q{trivializing operator} for operators introduced by Melrose and Piazza.

Let $A$ be a self-adjoint regular operator with compact resolvents.
We say that a self-adjoint operator $C$ of finite range 
is an \textit{$r$-trivializing operator} (or simply a trivializing operator)
for $A$ if the sum $A+C$ is invertible and the range of $C$ lies in the range of $\chirr(A)$.
\index{operators!trivializing}

Obviously, if $C$ is an $r$-trivializing operator for $A$, then 
\begin{equation}\label{eq:PAC}
	P = \chip(A+C)
\end{equation}
is an $r$-spectral section for $A$.
We say that a trivializing operator $C$ \textit{agrees} with a spectral section $P$ if \eqref{eq:PAC} holds.

These notions are generalized to the family case in a natural way.
We say that a norm continuous family $\C\colon X\to\Bsa$ 
is a \textit{trivializing family} for $\A\colon X\to\CRRsa$
if there is a continuous function $r\colon X\to\R_+$ such that 
$\C_x$ is an $r_x$-trivializing operator for $\A_x$ for every $\xX$.
We also call such a family $\C$ an \textit{$r$-trivializing} family for $\A$.
We say that a trivializing family $\C$ \textit{agrees} with a spectral section $P$ 
if $P_x = \chip(\A_x+\C_x)$ for every $\xX$.

\sub{Trivializing operators and spectral sections.}
Let $\T$ be the subspace of $^r\CRRsa\times\Bsa\times\R_+$ consisting of triples 
$(A,C,r)$ such that $C$ is an $r$-trivializing operator for $A$.

Recall that $\S$ denotes the subspace of $^r\CRRsa\times\Proj(H)\times\R_+$ 
consisting of triples $(A,P,r)$ such that $P$ is an $r$-spectral section for $A$.
There is the natural map
\begin{equation}\label{eq:DS}
	\tau\colon\T\to\S, \quad \tau(A,C,r) = (A,P,r), \text{ where } P = \chip(A+C).
\end{equation}
\index{3@Bundles and bundle maps!3tau@$\tau\colon\T\to\S$}

\uppskip
\begin{prop}\label{prop:DScont}
 The map \eqref{eq:DS} is continuous.
\end{prop}

\proof
The map $S\colon\T\to\Proj(H)$ taking $(A,C,r)$ to $\chill(A)$ is continuous on 
\[ \T_{\lambda} = \sett{(A,C,r)\in \T}{\pm\lambda\in\Res(A) \text{ and } r<\lambda }. \]
The positive spectral projection $P$ of $A+C$ can be written as an orthogonal sum 
\[ P = \chip(SAS+C) + \chilp(A) = \chip\br{SAS+C-(1-S)} + \chilp(A). \]
The operator $SAS+C-(1-S)$ is bounded, invertible, and depends continuously on $(A,C,r)\in\T_{\lambda}$.
Therefore, $P$ depends continuously on $(A,C,r)\in\T_{\lambda}$.
Since the open sets $\T_{\lambda}$ cover $\T$, the map \eqref{eq:DS} is continuous.
\endproof

\sub{Convexity.}
The fiber $\tau\inv(A,P,r)$ can be identified with a subset of $\Bsa$,
\begin{equation*}\label{eq:piC}
	\sett{C\in\Bsa}{\tau(A,C,r)=(A,P,r)}.
\end{equation*}

\begin{prop}\label{prop:Dconv}
 Under this identification, every fiber of $\tau$ is a convex subset of $\Bsa$.
\end{prop}

\proof
Let $C_0, C_1\in\tau\inv(A,P,r)$ and $C_t=(1-t)C_0+tC_1$, $t\in[0,1]$.
Then the range of $C_t$ lies in the range of $\one_{(-r,r)}(A)$.
Both $\A+C_0$ and $A+C_1$ commute with $P$, 
so the sum $A+C_t = (1-t)(A+C_0) + t(A+C_1)$ also commutes with $P$.
By the same reasoning, $A+C_t$ is strictly positive on the range of $P$ and strictly negative on the kernel of $P$.
Therefore, $C_t$ is an $r$-trivializing operator for $A$ which agrees with $P$,
that is, $(A,C_t,r)\in\tau\inv(A,P,r)$ for every $t\in[0,1]$.
\endproof

\begin{thm}\label{thm:DS}
 The bundle $\tau\colon\T\to\S$ is shrinkable
 (that is, fiber homotopy equivalent to the trivial bundle $\S\to\S$). 
 Moreover, a section 
\begin{equation}\label{eq:gamma}
	(A,P,r)\mapsto(A,\gamma(A,P,r),r) 
\end{equation}
 of $\tau$ can be chosen so
 that $\norm{\gamma(A,P,r)}<2r$ for every $(A,P,r)\in\S$.
\index{2@Maps!3gam@$\gamma\colon\S\to\Bsa$}
\end{thm}

\proof
If $\tau$ has a section, say \eqref{eq:gamma}, 
then the composition of this section with $\tau$ is fiber homotopic to the identity map $\Id_{\T}$
by Proposition \ref{prop:Dconv}, via the linear homotopy
$h_t(A,C,r) = (A,C_t,r)$, $C_t = (1-t)\cdot\gamma(\tau(A,C,r)) + t\cdot C$.
Therefore, we only need to prove the second statement of the theorem.

Let us fix a continuous function $\psi\colon\R\to[0,1]$ 
which is equal to $1$ on $(-\infty,0]$ and to $0$ on $[1,+\infty)$.
Our construction of the map $\gamma\colon\S\to\Bsa$ depends 
on the choice of such a function $\psi$.

Let $(A,P,r)\in\S$.
We take $Q=1-P$, $A^+ = PA P$, and $A^- = QA Q$.
The operators $r+A^+$ and $r-A^-$ are positive and invertible, 
and the sum $A^+ + A^-$ is orthogonal. 
The difference 
\begin{equation}\label{eq:C'}
	C' = (A^+ + A^-) - A = - PAQ - QAP
\end{equation}
is a self-adjoint operator of norm less than $r$
with the range of $C'$ lying in the range $\Ho$ of $\chirr(A)$.
Let 
\begin{equation}\label{eq:C''}
  C^+ = P\,\psi\br{r\inv A^+}P, \quad
   C^- = - Q\,\psi\br{-r\inv A^-}Q, \; \text{ and } \;
	 C'' = r\br{C^+ + C^-}.
\end{equation}
The range of $C''$ lies in $\Ho$ and $\norm{C''}\leq r$,
so the range of $C = C' + C''$ lies in $\Ho$ and $\norm{C}<2r$.

We define $\gamma$ by the formula $\gamma(A,P,r) = C = C' + C''$,
where $C'$ and $C''$ are defined by formulas \eqref{eq:C'} and \eqref{eq:C''}.

The function $\Psi(t) = t+\psi(t)$ is strictly positive on $(-1,+\infty)$. 
The sum $A+C$ can be written as
\begin{equation*}\label{eq:A'}
	A+C = (A^++A^-) + r\br{C^++C^-} = 
	 r\,P\,\Psi(r\inv A^+)P - r\,Q\,\Psi(-r\inv A^-)Q,
\end{equation*}
so it is invertible and $\chip(A+C) = P$.

It remains to show that $\gamma$ is continuous.
For $\lambda>0$, let
\[	\S_{\lambda} = \sett{(A,P,r)\in \S}{\pm\lambda\in\Res(A) \text{ and } r<\lambda }. \]
The projection $S = \chill(A)$ and the cut-off operator 
$B = SAS = f\inv\br{ S f(A) S } \in\B\sa(\H)$ 
depend continuously on $(A,P,r)\in \S_{\lambda}$.
The projections $S$ and $P$ commute, so $C' = S C' S = PBQ + QBP$, 
and thus the restriction of $C'$ to $\S_{\lambda}$ is continuous.
The restriction of $C''$ to $\S_{\lambda}$ can be written in a similar manner:
\[
 C'' = r\,S\brr{ P\psi(r\inv PBP)P - Q\psi(-r\inv QBQ)Q}S, 
\]
which provides continuity of $C''$ on $\S_{\lambda}$.
Open sets $\S_{\lambda}$ cover $\S$ when $\lambda$ runs $\R_+$, 
so both $C'$ and $C''$ are continuous on the whole $\S$.
Therefore, the sum $\gamma = C'+C''$ is also continuous.
This completes the proof of the first two statements of the theorem.
\endproof

\begin{thm}\label{thm:ss-tr-op}
  Let $X$ be an arbitrary topological space and $\A\colon X\to\CRRsa$ be a Riesz continuous map.
	Suppose that $\A$ admits an $r$-spectral section $P$.
	Then $\A$ admits an $r$-trivializing family $\C\colon X\to\Bsa$ that agrees with $P$
	and satisfies $\norm{\C_x}<2r_x$.	
\end{thm}

\proof
Let $\gamma\colon\S\to\Bsa$ be a map satisfying conditions of Theorem \ref{thm:DS}.
Then the formula $\C_x = \gamma(\A_x,P_x,r_x)$ determines an $r$-trivializing family for $\A$ 
satisfying conditions of the theorem.
\endproof

\begin{thm}\label{thm:ss-tr-op-c}
  Let $P$ be a spectral section for a Riesz continuous map $\A\colon X\to\CRRsa$.
	Then the space of all trivializing families for $\A$ that agree with $P$ 
	is convex and therefore contractible.		  
\end{thm}

\proof
Let $\C_{i}$ be an $r_{i}$-trivializing family for $\A$, $i=0,1$.
Then both $\C_0$ and $\C_1$ are $r$-trivializing families for $\A$, where $r=\max(r_0,r_1)$.
By Proposition \ref{prop:Dconv}, for every $t\in[0,1]$, $\C = (1-t)\C_0 + t\C_1$ 
is an $r$-trivializing family for $\A$ as well.
This completes the proof of the theorem.
\endproof

\section{$\Z_2$-graded case}\label{sec:Cl1}

Throughout this section $H = H^0\oplus H^1$ will be a $\Z_2$-graded Hilbert space.
Let $\sigma = \smatr{1 & 0 \\ 0 & -1}\in\B(H)$ be the symmetry defining the grading.
\index{5s@$\sigma$, grading}
We denote by $\CRRg$ the subset of $\CRRsa$ consisting of odd operators,
that is, operators anticommuting with $\sigma$.\index{1@Spaces!1R3o@$\Reg^1(H)$, odd self-adjoint regular operators}
Similarly, we denote by $\B^1(H)$ the subspace of $\B\sa(H)$ consisting of odd operators.\index{1@Spaces!1B2sao@$\B^1(H)$, odd self-adjoint operators}

\sub{Spectral sections for odd operators.}
The natural inclusion $\CRRg\hookto\CRRsa$ admits a spectral section $P\colon\CRRg\to\Proj(H)$.
Indeed, fix a continuous even function $\psi\colon\R\to\R$ supported on $[-r,r]$
which does not vanish at zero.
For every $A\in\CRRg$, the self-adjoint finite rank operator $C_A = \sigma\,\psi(A)$ 
is $r$-trivializing for the operator $A$, since $(A+C_A)^2 = A^2+\psi(A)^2$ is invertible.
Therefore, $P_A = \chip(A+C_A)$ is an $r$-spectral section for $A$.
Moreover, the maps $A\mapsto C_A$ and $A\mapsto P_A$ are continuous on $\CRRg$ 
and thus determine an $r$-trivializing family and an $r$-spectral section for the inclusion $\CRRg\hookto\CRRsa$.

It follows that, by a trivial reason, \textit{every} family of \textit{odd} self-adjoint operators with compact resolvents 
has a spectral section. 
Hence the notion of a spectral section is not very relevant for such operators.
Instead, one should consider spectral sections behaving well with respect to the grading.
Such a notion of a $\Cl(1)$ spectral section was introduced by Melrose and Piazza in \cite{MP2}.

\sub{$\Cl(1)$ spectral sections.}
A norm continuous family of projections $P\colon X\to\Proj(H)$ is called 
a $\Cl(1)$ spectral section for a family $\A\colon X\to\CRRg$ 
if $P$ is a spectral section for $\A$ and satisfies additionally the anticommutation property 
\begin{equation}\label{eq:P1}
	\sigma P\sigma\inv = 1-P.
\end{equation}
\index{spectral section!$\Cl(1)$}
Let $\Pg(H)$ denote the subspace of $\Proj(H)$ consisting of projections $P$ 
satisfying \eqref{eq:P1}.
Equivalently, a projection $P$ lies in $\Pg(H)$ if the symmetry $2P-1$ anticommutes with $\sigma$.

$\Pg(H)$ is naturally homeomorphic to the space of unitary operators $\U(H_0,H_1)$;
the corresponding homeomorphism 
\begin{equation}\label{eq:nu}
	\nu\colon\Pg(H)\to\U(H_0,H_1)
\end{equation}
takes $P\in\Pg(H)$ to $v\in\U(H_0,H_1)$ such that $2P-1 = \smatr{0 & v^* \\ v & 0}$.

\sub{Generalized $\Cl(1)$ spectral sections.}
Again, we consider two cases in parallel:
\begin{enumerate}[topsep=-1pt, itemsep=0pt, parsep=3pt, partopsep=0pt]
	\item $\ZZ = {^r\CRRsa}$, $\ZZ' = {^r\CRR(H^0,H^1)}$, and $I=\R_+$;
	\item $\ZZ=\Beusa$, $\ZZ'=\Beu(H^0,H^1)$, and $I=(0,1)$.
\end{enumerate}
\index{1@Spaces!2Z@$\ZZ$}
\index{1@Spaces!2Z@$\ZZ'$, $\Zg$}
Let $\Zg$ denote the subspace of $\ZZ$ consisting of odd operators.
The formula 
\begin{equation}\label{eq:Ahat}
	A\mapsto\hat{A} = \matr{0 & A^* \\ A & 0}
\end{equation}
defines a natural homeomorphism $\ZZ'\to\Zg$.
\index{5A@$\hat{A}$}

Both $\ZZ'$ and $\Zg$ are empty if one of $H^i$ is finite-dimensional
(recall that $H$ itself is infinite-dimensional). 
So we will always suppose that both $H^0$ and $H^1$ are infinite-dimensional.

We define a \textit{generalized $\Cl(1)$ spectral section} for $\hat{\A}\colon X\to\Zg$
as a generalized spectral section $P$ for $\hat{\A}$ satisfying \eqref{eq:P1}.
\index{spectral section!generalized $\Cl(1)$}
Equivalently, a norm continuous map $P\colon X\to\Pg(H)$ 
is a generalized $\Cl(1)$ spectral section for $\hat{\A}$ 
if $a_x$ is a compact deformation of the unitary $\nu(P_x)$ for every $\xX$,
where $a_x=A_x$ for $\ZZ=\Beusa$ and $a_x=f(A_x)$ for $\ZZ=\CRRsa$.

An element $a\in\Beu(H^0,H^1)$ is a compact deformation of a unitary 
if and only if the index of $a$ vanishes.
For $A\in\CRR(H^0,H^1)$ the indices of $A$ and f(A) coincide.
Therefore, the index of $A\in\ZZ'$ vanishes if and only if the operator $\hat{A}$ 
given by formula \eqref{eq:Ahat} admits a generalized $\Cl(1)$ spectral section.
We denote by $\bar{\ZZ}'$ the subspace of $\ZZ'$ consisting of operators with vanishing index
and by $\Zgg$ the subspace of $\Zg$ consisting of operators admitting a generalized $\Cl(1)$ spectral section.
\index{1@Spaces!2Z@$\bar{\ZZ}'$, $\Zgg$}

Since the index is locally constant on $\Beu(H^0,H^1)$ and thus also on $^r\CRR(H^0,H^1)$, 
$\bar{\ZZ}'$ is an open and closed subspace of $\ZZ'$
(in fact, it is a connected component of $\ZZ'$, but we do not use its connectedness).
Homeomorphism \eqref{eq:Ahat} takes $\bar{\ZZ}'$ to $\Zgg$,
so $\Zgg$ is an open and closed subspace of $\Zg$.

\begin{prop}\label{prop:Z1}
  For every $A\in\ZZ'$, the following four conditions are equivalent:
\begin{enumerate}[topsep=-3pt, itemsep=0pt, parsep=3pt, partopsep=0pt]
	\item The index of $A$ vanishes.
	\item The signature of the restriction of $\sigma$ to the kernel of $\hat{A}$ vanishes.
	\item $\hat{A}$ has a $\Cl(1)$ spectral section.
	\item $\hat{A}$ has a generalized $\Cl(1)$ spectral section.
\end{enumerate}
\end{prop}

\proof
($3\Rightarrow 4$) is trivial.

($1\Leftrightarrow 4$) is explained above.

($1\Leftrightarrow 2$) follows from the equality 
\begin{multline*}
	\sign\br{\restr{\sigma}{\Ker\hat{A}}}
  = \dim\br{H^0\cap\Ker\hat{A}} - \dim\br{H^1\cap\Ker\hat{A}} \\
  =  \dim(\Ker A) - \dim(\Ker A^*) = \ind(A).
\end{multline*}

($1\Rightarrow 3$)
If $\dim\Ker A = \dim\Ker A^*$, then there is a unitary $v\in\U(\Ker A, \Ker A^*)$.
Let $S\in\Proj(H)$ be the orthogonal projection onto the kernel of $\hat{A}$
and $P_0=\smatr{0 & v^* \\ v & 0}\in\Proj(\Ker\hat{A})$.
Then $P = P_0 S + \one_{(0,+\infty)}(\hat{A})$ is a $\Cl(1)$ spectral section for $\hat{A}$.
\endproof

\begin{prop}\label{prop:Z1r}
  For every $A\in\Zg$ and $r\in I$,
	the signatures of the restrictions of $\sigma$ to $\Ker(A)$ 
	and to the range of $\chirr(A)$ coincide.
\end{prop}

\proof
The range $V$ of $\chirr(A)$ can be decomposed into the orthogonal sum 
$V = V^-\oplus\Ker(A)\oplus V^+$ corresponding to the decomposition 
of the interval $(-r,\,r) = (-r,\,0)\cup\set{0}\cup(0,\,r)$.
Since $\sigma$ anticommutes with $A$, $\sigma$ takes $V^-$ to $V^+$ and vice versa.
Therefore, the signature of the restriction of $\sigma$ to $V^-\oplus V^+$ vanishes,
and so the signatures of $\restr{\sigma}{V}$ and $\restr{\sigma}{\Ker(A)}$ coincide.
\endproof

\sub{Homotopy Lifting Property.}
Let $\Gg$ denote the subspace of $\Zg\times\Pg(H)$ consisting of pairs $(A,P)$ 
such that $P$ is a generalized spectral section for $A$.
Recall that $\Zgg$ denotes the subspace of $\Zg$ 
consisting of operators admitting a generalized $\Cl(1)$ spectral section;
in other words, $\Zgg$ is the range of the map $\pig\colon\Gg\to\Zg$ taking $(A,P)$ to $A$.
Let 
\[ \UK(H^0) = \sett{u\in\U(H^0)}{u-1\in\K(H^0)} \] 
be the group of unitaries which are compact deformation of the identity.
\index{1@Spaces!1UK@$\UK(H)$}

\begin{thm}\label{thm:Omega1}
	The map $\pig\colon\Gg\to\Zg$ is a locally trivial principal $\UK(H^0)$-bundle.
	Both $\Gg\to\Zgg$ and $\Gg\to\Zg$ have the Homotopy Lifting Property for all spaces.
\end{thm}
\index{3@Bundles and bundle maps!2piG1@$\pig\colon\Gg\to\Zg$, $\Gg\to\Zgg$}

\proof 
We define the action $\mu$ of $\UK(H^0)$ on the product $\Zgg\times\Pg(H)$ by the formula
\[ \mu_u(A,P) = \br{A,\, \nu\inv(\nu(P)u)}, \quad u\in\UK(H^0), \]
where $\nu\colon\Pg(H)\to\U(H_0,H_1)$ is homeomorphism \eqref{eq:nu}.
For every $P,P'\in\Pg(H)$, 
$P-P'\in\K(H)$ is equivalent to $\nu(P)-\nu(P')\in\K(H^0,H^1)$,
which in turn is equivalent to $\nu(P')\inv\nu(P)-1 \in\K(H^0)$.
Therefore, $\Gg$ is fixed by the action of $\mu$, 
and $\mu$ acts transitively on the fibers of $\Gg\to\Zgg$.
Obviously, this action is free.
It follows that $\Gg\to\Zgg$ is a principal $\UK(H^0)$-bundle.

The next step of the proof is local triviality of $\Gg\to\Zgg$.
Since this bundle is principal, it is sufficient to show that 
it allows a local section over a neighbourhood of an arbitrary	 point $A_0\in\Zgg$.
Fix $r>0$ such that $\pm r\in\Res(A_0)$.
Let $S^{\circ}_r(A)$ and $S^+_r(A)$ be the projections defined by formulae \eqref{eq:Sr}.
They are continuous on the neighbourhood $V_r = \sett{A\in\Zgg}{\pm r\in\Res(A)}$ of $A_0$.
The projections $\frac{1\pm\sigma}{2}S^{\circ}_r(A)$ are also continuous on $V_r$,
so their ranges are locally trivial vector bundles over $V_r$.
Let $V\subset V_r$ be a neighbourhood of $A_0$ over which these two vector bundles are trivial;
denote their restrictions to $V$ by $E^+$ and $E^-$.
The ranks of $E^+$ and $E^-$ are equal to the dimensions of their fibers over $A_0$;
since $A_0\in\Zgg$, these ranks coincide.
Choose a unitary bundle isomorphism $v\colon E^+\to E^-$.
Let $T = \smatr{0 & v^* \\ v & 0}$ be the corresponding odd symmetry
and $P = (T+1)/2 \in\Proj(E^+\oplus E^-)$ the bundle projection. 
Then the formula $A\mapsto \br{A,\, P(A)S^{\circ}_r(A) + S^{+}_r(A)}$
defines a section of $\Gg$ over $V$.
This completes the proof of the second step.

The same reasoning as in the end of the proof of Theorem \ref{thm:Omega} shows that 
a locally trivial bundle $\Gg\to \Zgg$ has the Homotopy Lifting Property for all spaces.
Since $\Zgg$ is closed and open in $\Zg$, the same is true for $\Gg\to \Zg$.
This completes the proof of the theorem.
\endproof

\begin{thm}\label{thm:gss-inv1}
  Let $X$ be an arbitrary topological space and $\ZZ$ be either $^r\CRRsa$ or $\Beusa$.
	Let $\A\colon X\to \Zg$ be a continuous map.
	Then the following two conditions are equivalent:
\begin{enumerate}[topsep=-3pt, itemsep=0pt, parsep=3pt, partopsep=0pt]
	\item $\A$ has a generalized $\Cl(1)$ spectral section.
	\item $\A$ is homotopic, as a map from $X$ to $\Zg$, to a family of invertible operators.
\end{enumerate}
\end{thm}

\proof
The proof reproduces completely the proof of Theorem \ref{thm:gss-inv},
 with (generalized) spectral sections replaced by (generalized) $\Cl(1)$ spectral sections 
 and Theorem \ref{thm:Omega1} used instead of Theorem \ref{thm:Omega}. 
 The only additional care is needed for the choice of a compact operator $K$:
 it should commute with $\sigma$.
 Then $T'(x)$ is odd and the homotopy $h$ consists of odd operators.
\endproof

\sub{Fiber homotopy equivalence.}
From now on till the end of the section 
$\ZZ$ denotes the space $\CRRg$ equipped with the Riesz topology.
Let $\Sg$ denote the subspace of $\S$ consisting of triples $(A,P,r)$ such that 
$A$ is odd and $P$ is a $\Cl(1)$ spectral section for $A$ with a cut-off parameter $r$.
Restricting $\pis$ to $\Sg$, we obtain a fiber bundle $\pis\colon\Sg\to\Zgg$.\index{3@Bundles and bundle maps!2piS@$\pis\colon\Sg\to\Zgg$}

\begin{thm}\label{thm:ss-gs1}
  The bundle map $\iota\colon\Sg\to\Gg$ taking $(A,P,r)$ to $(A,P)$
	is a fiber homotopy equivalence over $\Zgg$.
\index{3@Bundles and bundle maps!3i@$\iota\colon\Sg\to\Gg$}
	Moreover, for every $\eps>0$, a fiber-homotopy inverse bundle map 
	$\varphi = \varphi_{\eps} \colon\Gg\to\Sg$ 
	can be chosen so that $\norm{Q-P} < \eps$ for every $(A,P)\in\Gg$ and $(A,Q) = \iota\circ\varphi(A,P)$.
	\index{3@Bundles and bundle maps!4phi1@$\varphi_{\eps}\colon\Gg\to\Sg$}
\end{thm}

\proof
We will show that the bundle map $\varphi$ from the proof of Theorem \ref{thm:ss-gs}
maps the subspace $\Gg$ of $\Gs$ to the subspace $\Sg$ of $\S$,
and that the restriction of $\varphi$ to $\Gg$ is fiber homotopy inverse to $\iota\colon\Sg\to\Gg$.

We use the designations from the proof of Theorem \ref{thm:ss-gs}.
It will be convenient to use the following convention:
if $B$ is a self-adjoint operator, then we write $\tilde{B}$ as an abbreviation for $2B-1$.
We will also use the \q{sign function} 
\[ \rho\colon\R\setminus\set{0}\to\set{-1,1}, \quad 
\rho = \one_{(0,\,+\infty)} - \one_{(-\infty,0)}. \]
\index{2@Maps!1rho@$\rho$, sign function}
Let $(A,P)$ be an arbitrary element of $\Gg$.
Equality \eqref{eq:Tr} can be written equivalently as 
\[ \tilde{T}_{r}(A,P) = (\Srp(A)-\Srm(A)) + \Sro(A)\,\tilde{P}\,\Sro(A). \]
Since $A$ and $\tilde{P}$ are odd, $\tilde{T}_{r}(A,P)$ is also odd.
It follows that $\tilde{T}(A,P)$ and $\tilde{Q}(A,P) = \rho\br{\tilde{T}(A,P)}$ are odd.
Therefore, $\varphi$ takes $\Gg$ to $\Sg$.

Formula \eqref{eq:PQt} can be written as 
$2[P,Q]_t-1 = \rho\br{(1-t)\tilde{P}+t\tilde{Q}}$.
Therefore, all the homotopies constructed in the proof of Lemma \ref{lem:PQR} preserve $\Pg(H)$,
and thus the restriction of $\varphi$ to $\Gg$ is fiber homotopy inverse to the restriction of $\iota$ to $\Sg$.
This completes the proof of the theorem.
\endproof

\begin{thm}\label{thm:ss1}
  Let $X$ be an arbitrary topological space and $\A\colon X\to\CRRg$ be a Riesz continuous map.
  Then the following three conditions are equivalent:
\begin{enumerate}[topsep=-3pt, itemsep=0pt, parsep=3pt, partopsep=0pt]
	  \item $\A$ has a $\Cl(1)$ spectral section.
	  \item $\A$ has a generalized $\Cl(1)$ spectral section.
	  \item $\A$ is homotopic, via a Riesz continuous homotopy $X\times[0,1]\to\CRRg$, 
	             to a family of invertible operators.
\end{enumerate}
	If $P$ is a generalized $\Cl(1)$ spectral section for $\A$ and $\eps>0$, then
	a $\Cl(1)$ spectral section $Q$ for $\A$ can be chosen so that
	$\norm{Q-P}_{\infty} < \eps$ and $Q$ is homotopic to $P$ as a generalized $\Cl(1)$ spectral section.
\end{thm}

\proof
($1\Rightarrow 2$) is trivial.

($2\Leftrightarrow 3$) follows from Theorem \ref{thm:gss-inv1}.

($2\Rightarrow 1$) and the last part of the theorem follow from Theorem \ref{thm:ss-gs1},
in the same manner as in the proof of Theorem \ref{thm:ss}.
\endproof

\sub{Trivializing operators.}
Recall that in the previous section we denoted by $\T$ 
the subspace of $^r\CRRsa\times\Bsa\times\R_+$ consisting of triples 
$(A,C,r)$ such that $C$ is an $r$-trivializing operator for $A$.
Let $\Tg$ be the subspace of $\T$ consisting of triples $(A,C,r)$ with $A$ and $C$ \textit{odd} operators.
In other words, $\Tg$ is the subspace of $^r\CRRg\times\Bg(H)\times\R_+$ 
consisting of triples $(A,C,r)$ such that $C$ is an $r$-trivializing operator for $A$.

If $(A,C,r)\in\Tg$, then $P = \chip(A+C)$ is a $\Cl(1)$ spectral section for $A$.
Therefore, the restriction of $\tau\colon\T\to\S$ to $\Tg$ defines a natural projection 
$\tau\colon\Tg\to\Sg$.
\index{3@Bundles and bundle maps!3tau@$\tau\colon\Tg\to\Sg$}

\begin{thm}\label{thm:DS1}
 The bundle $\tau\colon\Tg\to\Sg$ is shrinkable.
 Moreover, its section $(A,P,r)\mapsto(A,\gamma(A,P,r),r)$ can be chosen so that 
$\norm{\gamma(A,P,r)}<2r$ for every $(A,P,r)\in\S$.
\index{2@Maps!3gam1@$\gamma\colon\Sg\to\Bg(H)$}
\end{thm}

\proof
Let $\gamma\colon\S\to\Bsa$ be the map constructed in the proof of Theorem \ref{thm:DS}.
Then $\gamma$ takes $\Sg$ to $\Bg(H)$ and thus defines a section satisfying conditions of the theorem.
Indeed, if $A$ is odd and $P$ is a $\Cl(1)$ spectral section for $A$,
then the conjugation by $\sigma$ takes $A^+$ to $-A^-$ and vice versa.
Therefore, both $C'$ and $C''$ are odd, 
and thus $C = C' + C''$ is also odd for every $(A,P,r)\in\Sg$.
It follows that $\gamma$ satisfies the second condition of the theorem.

By the same reasoning as in the proof of Theorem \ref{thm:DS},
the composition of any section of $\Tg\to\Sg$ with $\tau$ is fiber homotopic to the identity map $\Id_{\Tg}$
via the linear homotopy, and thus the bundle $\Tg\to\Sg$ is shrinkable.
\endproof

\begin{thm}\label{thm:ss1-def}
  Let $X$ be an arbitrary topological space and $\A\colon X\to\CRRg$ be a Riesz continuous map.
	Suppose that $\A$ admits a $\Cl(1)$ spectral section $P$ with a cut-off function $r$. 
	Then $\A$ admits an odd $r$-trivializing family $\C\colon X\to\Bg(H)$ that agrees with $P$
	and satisfies $\norm{\C_x}<2r_x$.	
\end{thm}

\proof
Let $\gamma\colon\Sg\to\Bg(H)$ be a map satisfying conditions of Theorem \ref{thm:DS1}.
Then the formula $\C_x = \gamma(\A_x,P_x,r_x)$ determines an odd $r$-trivializing family for $\A$ 
satisfying conditions of the theorem.
\endproof

\begin{thm}\label{thm:ss1-def-c}
  Let $P$ be a $\Cl(1)$ spectral section for a Riesz continuous map $\A\colon X\to\CRRg$.
	Then the space of all odd trivializing families for $\A$ that agree with $P$ is convex and therefore contractible.		  
\end{thm}

\proof
This space is the intersection, inside the vector space $C(X,\Bsa)$, 
of the vector subspace $C(X,\Bg(H))$ with the subset of all trivializing families for $\A$.
The last subset is convex by Theorem \ref{thm:ss-tr-op-c},
so their intersection is also convex.
\endproof

\sub{Non-self-adjoint operators.}
Passing from odd self-adjoint operators to their chiral components, 
we obtain the following result.

\begin{thm}\label{thm:nsa-def}
Let $X$ be an arbitrary topological space and $\A\colon X\to\CRR(H)$ be a Riesz continuous map.
Then the following two conditions are equivalent:
\begin{enumerate}[topsep=-3pt, itemsep=0pt, parsep=3pt, partopsep=0pt]
	\item The bounded transform $f\circ\A$ is a compact deformation of 
	             a norm continuous family of unitaries. 
	\item $\A$ is Riesz homotopic to a family of invertible operators.	              
\end{enumerate}
If this is the case, then there is a norm continuous family $\C=(\C_x)$ 
of finite rank operators such that $\A_x+\C_x$ is invertible for every $x\in X$, 
the range of $\C_x$ lies in the range of $\one_{[0,\,r_x)}(\A_x\A_x^*)$,
and the kernel of $\C_x$ contains the range of $\one_{[r_x,+\infty)}(\A_x^*\A_x)$
for some continuous function $r\colon X\to\R_+$ .
\end{thm}

\proof
The first part of the theorem is equivalent to the part $(2\Leftrightarrow 3)$ of Theorem \ref{thm:ss1} 
applied to the family  $\hat{\A} = \smatr{0 & \A^* \\ \A & 0}$ 
of regular odd self-adjoint operators with compact resolvents. 

Suppose now that $\A$ is Riesz homotopic to a family of invertible operators. 
Then, by Theorem \ref{thm:ss1}, $\hat{\A}$ has a $\Cl(1)$ spectral section;
let $R\colon X\to\R_+$ be its cut-off function.
Applying Theorem \ref{thm:ss1-def} to $\hat{\A}$, we get a norm continuous family 
$\hat{\C} = \smatr{0 & \C^* \\ \C & 0}$ of odd self-adjoint operators 
such that $\hat{\A}+\hat{\C}$ is invertible 
and the range of $\hat{\C}$ lies in the range of $\one_{(-R,\,R)}(\hat{\A})$
(we omit the subscript $x$ for brevity). 
Since
\[ \one_{(-R,\,R)}(\hat{\A}) \,=\, \one_{[0,\,R^2)}(\hat{\A}^2)
   = \one_{[0,\,R^2)}\br{\A^*\A\oplus\A\A^*}, \]
the family $\C$ and the function $r=\sqrt{R}$ satisfy conditions of the theorem.
\endproof

\addtocontents{toc}{\vspace{-5pt}}

\part{Special cases}\label{part:appl}

\addtocontents{toc}{\vspace{-12pt}}

In this part we present a number of special cases of general results obtained in the first part of the paper.

\section{Relatively compact deformations}\label{sec:appl-rc}

Let $H$ and $H'$ be Hilbert spaces.

\sub{Deformations of a single operator.}
We give here only two of possible examples.
Clearly, one can write a $\Z_2$-analog of Theorem \ref{thm:rel-comp},
using Theorems \ref{thm:ss1} and \ref{thm:ss1-def} instead of Theorems \ref{thm:ss} and \ref{thm:ss-tr-op};
we omit it since it is quite straightforward. 

\begin{thm}\label{thm:rel-comp-nsa}
For a regular operator $B\colon H\to H'$ with compact resolvents, let 
 \[ \ZZ_B = \sett{A\in\CRR(H,H')}{f(A)-f(B)\in\K(H,H')}. \] 
Then there are Riesz-to-norm continuous maps $\alpha\colon\ZZ_B\to\B(H,H')$ 
and $r\colon \ZZ_B\to\R_+$ such that $A+\alpha(A)$ is invertible, 
the range of $\alpha(A)$ lies in the range of $\one_{[0,\,r_x)}(AA^*)$,
and the kernel of $\alpha(A)$ contains the range of $\one_{[r_x,+\infty)}(A^*A)$
for every $A\in\ZZ_B$.
\end{thm}

\proof
The bounded transform $f(B)$ is a compact deformation of some unitary $u\in\U(H,H')$.
Therefore, the composition of the bounded transform  with the inclusion $\ZZ_B\hookto\CRR(H,H')$ 
is a compact deformation of a norm continuous (even constant!) family of unitaries
$\ZZ_B\ni A\mapsto u$.
It remains to apply Theorem \ref{thm:nsa-def}.
\endproof

\begin{thm}\label{thm:rel-comp}
For a regular self-adjoint operator $B\in\CRRsa$ with compact resolvents, let 
\begin{equation}\label{eq:ZB}
	\ZZ\sa_B = \ZZ_B\cap\CRRsa = \sett{A\in\CRRsa}{f(A)-f(B)\in\K(H)} \subset\CRR\sa(H)
\end{equation}
be the subspace of $\CRRsa$ equipped with the Riesz topology. 
Let $P = \chip(B)$ be the positive spectral projection of $B$. 
Then inclusion \eqref{eq:ZB} admits both a spectral section and a trivializing family. 
Moreover, for every $\eps>0$ a spectral section $Q$ can be chosen so that
$\norm{Q(A)-P} < \eps$ for every $A\in\ZZ\sa_B$.
\end{thm}

\proof
The constant map $A\mapsto P$ is a generalized spectral section for inclusion \eqref{eq:ZB}.
It remains to apply Theorems \ref{thm:ss} and \ref{thm:ss-tr-op}.
\endproof

\sub{Essentially self-adjoint operators.}
Recall that a bounded operator $a$ is called essentially self-adjoint
if $a-a^*$ is a compact operator.
\index{operators!bounded!essentially self-adjoint}

\begin{thm}\label{thm:ess-sa}
Let $X$ be the subspace of $\CRR(H)$ consisting of operators $A$ 
whose bounded transform is an essentially self-adjoint operator.
Then there are Riesz-to-norm continuous maps $\alpha\colon X\to\B(H)$ 
and $r\colon X\to\R_+$ such that $A+\alpha(A)$ is invertible, 
the range of $\alpha(A)$ lies in the range of $\one_{[0,\,r_x)}(AA^*)$,
and the kernel of $\alpha(A)$ contains the range of $\one_{[r_x,+\infty)}(A^*A)$
for every $A\in X$.
\end{thm}

\proof
Let $A\in X$ and $a=f(A)$. 
Then $b = (a+a^*)/2$ is a self-adjoint operator of norm $\leq 1$ 
and $u = u_A = b+i\sqrt{1-b^2}$ is a unitary.
Moreover, both $b$ and $u$ are compact deformations of $a$.
The map $X\to\U(H)$ taking $A$ to $u_A$ is Riesz-to-norm continuous.
Therefore, the composition of the bounded transform with the embedding $X\hookto\CRR(H)$
is a compact deformation of a norm continuous family of unitaries.
It remains to apply Theorem \ref{thm:nsa-def}.
\endproof

\sub{Essentially odd operators.}
Let $H$ be a $\Z_2$-graded Hilbert space, with the grading given by the symmetry $\sigma$.
A bounded operator $a\in\B(H)$ is called essentially odd if 
$\sigma a + a\sigma$ is a compact operator.
\index{operators!bounded!essentially odd}

\begin{thm}\label{thm:ess-odd}
Let $X$ be the subspace of ${^r}\CRR\sa(H)$ consisting of operators $A$ 
whose bounded transform is an essentially odd operator.
Then the natural embedding $X\hookto{^r}\CRR\sa(H)$ 
admits both a spectral section $X\to\Proj(H)$ and a trivializing family $X\to\Bsa$.
\end{thm}

\proof
Let $A\in X$ and $a=f(A)$. 
Then $b = (a - \sigma a\sigma)/2$ is an odd self-adjoint operator of norm $\leq 1$, 
$u = b + \sigma\sqrt{1-b^2}$ is a symmetry, and $P=P_A=(u+1)/2$ is a projection.
Moreover, both $b$ and $u$ are compact deformations of $a$.
The map $X\to\Proj(H)$ taking $A$ to $P_A$ is Riesz-to-norm continuous.
Therefore, this map is a generalized spectral section for the embedding $X\hookto{^r}\CRR\sa(H)$.
It remains to apply Theorems \ref{thm:ss} and \ref{thm:ss-tr-op}.
\endproof

\section{Pseudodifferential operators}\label{sec:PDO}

We show here several examples of applications of our results 
to pseudodifferential operators over closed manifolds. 
Again, we omit the $\Z_2$-graded case here. 

Let $M$ be a closed smooth manifold equipped with a smooth positive measure, 
and let $E$, $E'$ be smooth Hermitian bundles over $M$.
Let $\Psi_d(E,E')$ denote the space of pseudodifferential operators 
of order $d\geq 0$ acting from sections of $E$ to sections of $E'$.
\index{1@Spaces!3Psi@$\Psi_d(E,E')$, pseudodifferential operators of order $d$}
We equip it with the topology induced by the inclusion 
\[ \Psi_d(E,E')\hookto\B(L^2_d(M;E),L^2(M;E'))\times\B(L^2_d(M;E'),L^2(M;E)) \]
taking a pseudodifferential operator $A$ to the pair $(A,A^t)$,
where $A^t$ is the operator formally adjoint to $A$
and $L^2_d(M;E)$ is the order $d$ Sobolev space of sections of $E$.
\index{1@Spaces!3Hd@$L^2_d(E)$, Sobolev space}
The natural inclusion of the subspace 
$\Psi\el_d(E,E')\subset\Psi_d(E,E')$ of elliptic operators to $\Reg(L^2(M;E),L^2(M;E'))$
is Riesz continuous (see, for example, \cite[Proposition 2.2]{L04}).
\index{1@Spaces!3Psi2@$\Psi\el_d(E,E')$, elliptic pseudodifferential operators}

\begin{thm}\label{thm:PDO-nsa}
  Let $X$ be a topological space and $\A\colon X\to\Psi\el_d(E,E')$ 
	be a continuous family of elliptic operators of order $d\geq 1$.
	Suppose that $\A$ is homotopic to a family of invertible operators.
	Then there is a norm continuous family $\alpha=(\alpha_{x,B})$ 
	of smoothing finite rank operators
	parametrized by pairs $(x,B)\in X\times\Psi_{d-1}(E,E')$
	such that $\A_x+B+\alpha_{x,B}$ is invertible for every $\xX$ and every $B\in\Psi_{d-1}(E,E')$.
\end{thm}

\begin{rem}
In this theorem, $\A$ is homotopic to a family of invertible operators, in particular, 
in each of the following cases:
\begin{enumerate}[topsep=-3pt, itemsep=0pt, parsep=3pt, partopsep=0pt]
	\item $X$ is contractible.
	\item $X$ is compact and $\ind(\A)=0\in\K^0(X)$.
  \item The kernel and cokernel of $\A_x$ have locally constant ranks,
	             and the corresponding vector bundles $\Ker(\A)$ and $\Coker(\A)$ over $X$
							 are isomorphic.
	\end{enumerate}
\end{rem}

\proof
Let $Y=X\times\Psi_{d-1}(E,E')$, H=$L^2(M;E)$, and $H'=L^2(M;E')$.
The map $\tilde{\A}\colon Y\to\Psi\el_d(E,E')$ 
taking $(x,B)$ to $\A_x+B$ is continuous.
Therefore, the composed map 
$$\tilde{\A}\colon Y\to\Psi\el_d(E,E')\hookto\CRR(H,H')$$
is Riesz continuous.
It is Riesz homotopic to the map $Y\to\CRR(H,H')$ taking $(x,B)$ to $\A_x$
via the homotopy $(x,B,t)\mapsto \A_x+tB$.
Since $\A$ is homotopic to a family of invertible operators,
the same is true for $\tilde{\A}$.

By Theorem \ref{thm:nsa-def}, there is a continuous map $\alpha\colon Y\to\B(H,H')$
such that $\tilde{\A}+\alpha$ is an invertible family, 
the range of $\alpha_{x,B}$ lies in the range $V$ of $\one_{[0,\,r)}(\tilde{\A}_{x,B}\tilde{\A}_{x,B}^*)$,
and the orthogonal complement of the kernel of $\alpha_{x,B}$ 
lies in the range $V'$ of $\one_{[0,\,r)}(\tilde{\A}_{x,B}^*\tilde{\A}_{x,B})$
for some $r=r(x,B)$.
Since $\tilde{\A}_{x,B}$ is an elliptic operator of positive order,
both $V$ and $V'$ are spanned by a finite number of $C^{\infty}$-sections.
Therefore, $\alpha_{x,B}$ is a smoothing operator of finite rank.
This completes the proof of the theorem.
\endproof

\sub{Self-adjoint case.}
For $M$ and $E$ as above, 
let $\Psi\sa_d(E)$ denote the subspace of $\Psi_d(E)$ consisting of symmetric operators,
and let $\Psi\elsa_d(E) = \Psi\el_d(E)\cap\Psi\sa_d(E)$.

\begin{thm}\label{thm:PDO-sa}
  Let $X$ be a topological space and $\A\colon X\to\Psi\elsa_d(E,E')$ 
	be a continuous family of symmetric elliptic operators of order $d\geq 1$.
	Suppose that $\A$ is homotopic to a family of invertible operators.
  Then the map 
	\[ \tilde{\A}\colon Y = X\times\Psi\sa_{d-1}(E)\to\Psi\elsa_d(E), \quad	(x,B)\mapsto\A_x+B \]
	admits both a spectral section $P\colon Y\to\Psi\sa_0(E)$ 
	and a trivializing family with smoothing trivializing operators.
\end{thm}

\proof
The proof is completely similar to the proof of Theorem \ref{thm:PDO-nsa};
one only needs to use Theorems \ref{thm:ss} and \ref{thm:ss-tr-op} instead of Theorem \ref{thm:nsa-def}.
\endproof

\begin{rem}
In this theorem, $\A$ is homotopic to a family of invertible operators, in particular, 
in each of the following cases:
\begin{enumerate}[topsep=-3pt, itemsep=0pt, parsep=3pt, partopsep=0pt]
	\item $X$ is contractible.
	\item $X$ is compact and $\ind(\A)=0\in\K^1(X)$.
  \item The kernel of $\A_x$ has locally constant rank.
\end{enumerate}
\end{rem}

\section{Cobordism theorems}\label{sec:cob}

Let $N$ be a smooth compact Riemannian manifold with non-empty boundary $M = \pN$
and $F$, $F'$ be smooth Hermitian vector bundles over $M$.
Denote by $E$ and $E'$ the restrictions of $F$ and $F'$ to $M$.

\sub{Calder\'on projection.}
Let $D$ be a first order elliptic differential operator over $N$ 
acting from sections of $F$ to sections of $F'$.
The \textit{space of Cauchy data} of $D$ is the closure in $H = L^2(M;E)$ 
of the subspace consisting of restrictions to $M$ of all smooth solutions of the equation $Du=0$.
\index{space of Cauchy data}
The orthogonal projection $Q=Q(D)\in\Proj(H)$ onto the Cauchy data space   
is called the (orthogonal) \textit{Calder\'on projection} of $D$;
it is a pseudodifferential operator of zeroth order.
\index{Calder\'on projection}

At the points of the boundary the operator $D$ can be written as 
\begin{equation}\label{eq:DJA}
	D = -iJ(\p_z + A), 
\end{equation}
where $z$ is the normal coordinate, 
$J$ is the conormal symbol of $D$ (a bundle isomorphism from $E$ to $E'$
),
and $A$ is a first order elliptic differential operator over $M$ 
acting on sections of $E$.
Such an operator $A$ is called \textit{the tangential operator of $D$ along the boundary},
or simply \textit{the boundary operator} of $D$.
\index{tangential operator}
\index{boundary operator}

Suppose that the principal symbol of $A$ is self-adjoint, that is $A-A^t$ is a bundle endomorphism
(here $A^t$ denotes the operator formally adjoint to $A$).
Then the Calder\'on projection $Q(D)$ has the same principal symbol as the positive spectral projection $\chip(A+A^t)$.
In other words, $Q(D)$ is a generalized spectral section for the symmetrized tangential operator $\At = (A+A^t)/2$
\index{symmetrized tangential operator},
or for any other symmetric operator on $E$ with the same principal symbol as $A$.

\begin{rem}\label{rem:DJA}
	Near the boundary, any first order operator $D$ can be written in the form 
\begin{equation}\label{eq:DJAz}
	D = -iJ(z)(\p_z + A(z)) 
\end{equation}
	where $(J(z))$ and $(A(z))$ are one-parameter families 
	of bundle endomorphisms and of differential (tangential) operators respectively,
	which depend on the parameter $z$.
	Expression \eqref{eq:DJA} should be understood as the restriction of \eqref{eq:DJAz} to the boundary $z=0$.
	In other words, $J$ and $A$ in \eqref{eq:DJA} are just the boundary values of $J(z)$ and $A(z)$ at $z=0$.

	In the studies of boundary value problems, 
	the operator $D$ is often required to have a product form near the boundary,
	which means that both $J(z)$ and $A(z)$ are actually independent of $z$.
	Since the homotopy methods that we use in this paper are insensitive to continuous perturbations, 
	we do not need this requirement.
	We do not require the product form of the operator near the boundary,
	neither we require the product form of the metric.
\end{rem}

\upskip

\sub{General cobordism theorem.}
The following result is the most general form of a cobordism theorem based on Calder\'on projections.
Below we apply Theorem \ref{thm:cob} to Dirac type operators, 
for which continuous dependence of the Calder\'on projection on parameter is known.

\begin{thm}\label{thm:cob}
Let $N$, $F$, and $F'$ be as above, $X$ be an arbitrary topological space, 
$\D=(\D_x)_{\xX}$ be a family of first order elliptic differential operators
$\D_x\colon C^{\infty}(N;F)\to C^{\infty}(N;F')$,
and $\A_x$ be the boundary operator of $\D_x$. 
Suppose that the principal symbol of each $\A_x$ is self-adjoint 
and that the Calder\'on projection $Q(\D_x)$ depends continuously on $x$.
Let $\AAt = (\At_x)$ be a continuous family 
of first order symmetric operators over $M$, with $\At_x$ having the same principal symbols as $\A_x$.
Then $\AAt$ admits both a spectral section 
and a trivializing family, with smoothing trivializing operators.
\end{thm}

\proof
The Calder\'on projection $Q_x = Q(\D_x)$ is a generalized spectral section for $\AAt_x$.
The norm continuous family $(Q_x)$ of projections is a generalized spectral section for a Riesz continuous family $\AAt$.
It remains to apply Theorems \ref{thm:ss} and \ref{thm:ss-tr-op}.
Since the range of every trivializing operator lies in the span of a finite number 
of $C^{\infty}$-sections (namely, eigenfunctions of $\AAt_x$),
all trivializing operators are smoothing. 
\endproof

\sub{Continuity of Calder\'on projections.}
There are different criteria of continuity of the family of Calder\'on projections
for different classes of operators.
One of such criteria is \cite[Corollary 7.4]{BBLZ}. 
However, this criteria uses a specific metric on the space of operators (\q{strong metric} of \cite[Definition 7.1]{BBLZ}), 
which makes it is difficult for use in applications.
We prefer to use another criteria, 
which is proven by Booss-Bavnbek, Deng, Zhou, and Zhu in \cite{BDZZ}. 
Unfortunately, loc. cit. does not contain the exact statement we need, so we state it here explicitly.

Recall that a first order operator $D$ is said to have weak inner Unique Continuation Property (weak inner UCP) 
\index{weak inner UCP}
if the only inner solution (that is, a solution vanishing on the boundary)
\index{inner solution}
of the equation $Du=0$ is the trivial solution $u=0$.

\begin{prop}[\cite{BDZZ}]\label{prop:Cald-cont}
Let $N$, $F$, and $F'$ be as above, $X$ be a topological space, 
and $\D=(\D_x)_{\xX}$ be a family of first order elliptic differential operators
$\D_x\colon C^{\infty}(N;F)\to C^{\infty}(N;F')$.
Suppose that the following two conditions hold:
\begin{enumerate}[topsep=-3pt, itemsep=0pt, parsep=3pt, partopsep=0pt]
	\item The coefficients of the operators $\D_x$ and $\D^t_x$, 
	together with the first derivatives of the conormal symbol of $\D_x$ in the tangential directions, 
depend continuously on $x$ in every local chart of $N$.
	\item All the operators $\D_x$ and $\D^t_x$ have weak inner UCP
(or, more generally, the dimensions of the spaces of inner solutions of $\D_x u=0$ and $\D_x^t u=0$ do not depend on $\xX$).
\end{enumerate}
Then the family of Calder\'on projections $x\mapsto Q(\D_x)\in \Proj(L^2(M;E))$ is norm continuous.
\end{prop}

\proof
The main theorem of \cite{BDZZ}, Theorem 1.2, is stated for operators of order $d$ in the following form:
if its conditions (i) and (ii) hold for \textit{every} $s\geq d/2$, 
then the family of orthogonal Calder\'on projections $(Q_x)$ is continuous in the operator norms 
of \textit{all} the Sobolev spaces $H^r(M;E)$, $r\in\R$.
However, their reasoning actually proves a stronger result:
if conditions (i) and (ii) hold for \textit{some} $s\geq d/2$, 
then $(Q_x)$ is continuous in the operator norms of $H^r(M;E)$ for every $r\in[-s,s]$.

We apply this stronger result to the family $(\D_x)$ with $d=1$, $s=1/2$, and $r=0$.
Condition (ii) coincides with condition (2) of our proposition.
The operator $\D_x$ depends continuously on $x$ in the operator norm $\norm{\cdot}_{1,\,0}$,
that is, the norm on the space of bounded operators $H^1(N;F)\to L^2(N;F')$.
The adjoint operator $\D^t_x$ has the analogous property.
The natural inclusion 
\[ C^1(\Hom(M;E),\Hom(M;E'))\hookto\B(H^{1/2}(M;E),H^{1/2}(M;E')) \]
is continuous, so the conormal symbol $\sigma_x(n)$ of $\D_x$ depends continuously on $x$ in the operator norm $\norm{\cdot}_{1/2,\,1/2}$.
Therefore, condition (i) of \cite[Theorem 1.2]{BDZZ} holds in our case for $s=1/2$.
It follows that $(Q_x)$ is a continuous family of projections in the operator norm of $L^2(M;E)$.
\endproof

\comment{============================================
\begin{rem*}
In the statement of the proposition, weak inner UCP can be replaced by a weaker condition, 
namely that the dimensions of the spaces of inner solutions of $\D_x u=0$ and $\D_x^t u=0$ do not depend on $\xX$.
But Proposition \ref{prop:Cald-cont} is sufficient for our purposes.
\end{rem*}
====================================================}

\sub{Unique Continuation Properties.}
Recall that an operator $D$ over a connected manifold $N$ is said to have weak Unique Continuation Property (weak UCP)
if any solution of the equation $Du=0$ which vanishes on an open subset of $N$ vanishes on the whole $N$.
\index{weak UCP}

Let $N$ be a smooth connected Riemannian manifold, not necessarily compact,
and $F$, $F'$ be smooth Hermitian vector bundles over $N$.
Denote by $\W(F,F')$ the set of first order elliptic differential operators 
$D\colon C^{\infty}(N;F)\to C^{\infty}(N;F')$
whose principal symbol $d$ satisfies the following condition:
\begin{equation}\label{eq:Wdef} 
\begin{split}
	& \text{The fiber endomorphism } i d(\xi)\inv d(\eta)\in\End(F_x) \text{ is self-adjoint} \\
	& \text{for every pair of orthogonal cotangent vectors } \xi,\eta\in T_x^*M, \; x\in M.
\end{split}
\end{equation}
\index{1@Spaces!3@$\W(F,F')$}

\begin{prop}\label{prop:UCP}
Every $D\in\W(F,F')$ has weak Unique Continuation Property.
\end{prop}

\proof
We follow the line of the proof of weak UCP for perturbed Dirac type operators 
in \cite[Theorem 1.33]{BBB}, but write it in more detail.

We can suppose without loss of generality that $N$ has no boundary 
(otherwise replace $N$ by $N\setminus\pN$).
Suppose that a nontrivial solution $u$ of $Du=0$ vanishes on a non-empty open subset of $N$.
Let $V$ be the union of all open subsets of $N$ where $u$ vanishes, 
and let $N' = \supp(u)\neq\emptyset$ be the complement of $V$ in $N$.

1. We claim that there is a point $p\in V$ such that the injectivity radius of $p$
is greater than the distance from $p$ to $N'$, $\inj(p) > \dist(p,N')$.
To show this, choose $x\in N'\cap\overline{V}$, and let $r=\inj(x)$.
Since the injectivity radius function is lower-semicontinuous,
there is a $\delta\in(0,r/4)$ such that the injectivity radius is greater than $r/2$ 
for all points of the open ball $B_{\delta}(x) = \sett{y\in N}{\dist(x,y)<\delta}$.
Since $x$ lies in the boundary of $V$, the intersection $V\cap B_{\delta}(x)$ is non-empty;
let $p$ be a point in this intersection.
Then $\dist(p,N')\leq\dist(p,x)<r/4$ and $\inj(p)>r/2$.

2. Let $p\in V$ be such a point that $r = \inj(p) > \dist(p,N') = d$. 
Then the open ball $B_d(p)$ is contained in $V$
and the larger open ball $B_r(p)$ can be equipped with (geodesical) spherical coordinates.
It follows from \cite[Lemmata 5 and 6]{BB} that 
$u$ vanishes on some intermediate ball $B_R(p)$ with $d<R<r$, that is $B_R(p)\subset V$.
On the other hand, the radius of $B_R(p)$ is greater than the distance from $p$ to $N'$,
so $B_R(p)$ intersects $N'$.
This contradiction shows that $D$ satisfies weak UCP,
which completes the proof of the proposition.
\endproof

\sub{Dirac type operators.}
Recall that a first order operator $D$ with the principal symbol $d$ 
is called a Dirac type operator if $d(\xi)^*d(\xi) = \norm{\xi}^2\cdot\Id = d(\xi)d(\xi)^*$ 
for every $\xi\in T^*M$.
\index{Dirac type operator}
Every Dirac type operator acting from sections of $F$ to sections of $F'$ is an element of $\W(F,F')$,
but not vice versa.

Proposition \ref{prop:UCP} can be used to prove a cobordism theorem for operators of the class $\W(F,F')$.
However, such a result is not really more general than a cobordism theorem for Dirac type operators.
Indeed, composing $D$ with a bundle automorphism of $F'$ 
does not affect the boundary operator of $D$.
The following proposition shows that every operator $D\in\W(F,F')$
can be obtain from a Dirac type operator by such a composition.

\begin{prop}\label{prop:W-Dirac}
Let $D$ be a first order operator acting from sections of $F$ to sections of $F'$.
Then the following two conditions are equivalent:
\begin{enumerate}[topsep=-2pt, itemsep=0pt, parsep=3pt, partopsep=0pt]
	\item $D\in\W(F,F')$,
	\item $D = I\cdot D'$, where $I$ is a bundle automorphism of $F'$ and $D'$ is a Dirac type operator.
\end{enumerate}
\end{prop}

\proof
Condition \eqref{eq:Wdef} can be equivalently written as follows:
\begin{equation}\label{eq:Wdef2} 
d(\xi)d(\eta)^* +d(\eta)d(\xi)^* = 0  \quad \text{for every } \xi\bot\eta\in T_x^*N, \; x\in N.
\end{equation}
Left hand side of \eqref{eq:Wdef2} is an $\End(F'_x)$-valued symmetric bilinear form on $T_x^*N$,
\[ \delta(\xi,\eta) = d(\xi)d(\eta)^* +d(\eta)d(\xi)^*. \]
($2\Rightarrow 1$)
If $D = I\cdot D'$, then $\delta(\xi,\eta) = 2\bra{\xi,\eta}I\,I^*$ satisfies \eqref{eq:Wdef2},
so $D\in\W(F,F')$.

($1\Rightarrow 2$)
Let $D\in\W(F,F')$.
For arbitrary non-zero $\xi,\eta\in T_x^*F$,
we write $\xi = \xi'+t\eta$ with $\xi'$ orthogonal to $\eta$. 
This gives
\begin{equation}\label{eq:xi-eta}
	\delta(\xi,\eta) = t\delta(\eta,\eta) = \bra{\xi,\eta}S(\eta),
\end{equation}
where 
$$S(\eta) = \delta(\eta,\eta)/\norm{\eta}^2 = 2d(\eta)d(\eta)^* / \norm{\eta}^2 \in\End(F'_x)$$ 
is a homogenous function on $S_x^*N\setminus\set{0}$ of degree 0.
Since $\delta(\xi,\eta) = \delta(\eta,\xi)$, identity \eqref{eq:xi-eta} impies 
$S(\xi) = S(\eta)$ for every pair of non-orthogonal vectors $\xi,\eta\in S_x^*N$.
For every non-zero $\xi,\eta\in T_x^*N$ there is a third vector $\zeta\in S_x^*N$
which is non-orthogonal to both $\xi$ and $\eta$, so that $S(\xi) = S(\zeta) = S(\eta)$. 
Therefore, $S(\eta)$ is independent of $\eta\in S_x^*N$ and depends only on $x$, $S(\eta)=S_x$.	
Moreover, $S_x$ is positive for every $x\in N$.
Let $I_x$ be the positive square root of $S_x/2$.
Then $I$ is a smooth bundle automorphism of $F'$. 
The equality $d(\xi)d(\xi)^* = \norm{\xi}^2 I_x I_x^*$ implies
$(I_x\inv d(\xi))(I_x\inv d(\xi))^* = \norm{\xi}^2$ for every $\xi\in T_x^*N$.
Since $I_x\inv d(\xi)$ is invertible for $\xi\neq 0$, this implies the second identity 
$(I_x\inv d(\xi))^*(I_x\inv d(\xi)) = \norm{\xi}^2$.
Therefore, $I\inv D$ is a Dirac type operator.
\endproof

\begin{prop}\label{prop:wiUCP}
Let $N$ be a smooth connected Riemannian manifold with non-empty boundary,
and let $D$ be a Dirac type operator over $N$.
Then $D$ has weak inner Unique Continuation Property.
\end{prop}

\proof
The operator $D$ admits an extension across the boundary,
that is, $D$ is the restriction to $N$ 
of some Dirac type operator $\Dt$ over $\Nt$, 
where $\Nt$ is a smooth Riemannian manifold without boundary 
containing $N$ as a smooth submanifold of codimension zero.
Indeed, let $D$ act from sections of $F$ to sections of $F'$.
The symbol $d$ of $D$ determines the structure of a Clifford module over $T^*N\oplus\R$ on $F\oplus F'$,
with the cotangent vector $\xi$ acting as $\hat{d}(\xi) = \smatr{0 & d(\xi)^* \\ d(\xi) & 0}$ 
and the unit vector in the additional $\R$-direction acting as $1_F\oplus(-1_{F'})$.
Such a structure can be smoothly extended across the boundary of $N$,
to an external collar neighbourhood of $\p N$ equipped with a metric in a compatible way.
Such an extension gives rise to the symbol $\tilde{d}$ of a Dirac type operator, 
which extends $d$ and acts from sections of $\Ft$ to sections of $\Ft'$,
where $\Ft$ and $\Ft'$ are smooth Hermitian vector bundles over $\Nt$ extending $F$ and $F'$ respectively.
Given an extension of the symbol, 
the whole operator $D$ can be smoothly extended to $\Nt$ using a partition of unity.

Let $u\in L^2(N;F)$ be in the kernel of $D$ 
and let $\tilde{u}\in L^2(\Nt;\Ft)$ be the extension of $u$ to $\Nt\setminus N$ by zero.
Suppose that $u$ vanishes on the boundary $\p N$.
Then the Green formula for $D$ implies that $\tilde{u}$ is a weak solution of $\Dt$. 
Since $\Dt$ is elliptic, $\tilde{u}$ is smooth.
By Proposition \ref{prop:UCP}, the operator $\Dt$ has weak UCP.
Since $\Dt\tilde{u} = 0$ and $\tilde{u}$ vanishes on the open subset $\Nt\setminus N$, 
we get $\tilde{u}\equiv 0$ and thus $u = \restr{\tilde{u}}{N}\equiv 0$.
Therefore, $D$ has weak inner UCP.
\endproof

\sub{Cobordism theorem for Dirac type operators.}
The boundary operator of a Dirac type operator is again a Dirac type operator;
moreover, it has a self-adjoint principal symbol.
Combining Theorem \ref{thm:cob} with Propositions \ref{prop:Cald-cont} and \ref{prop:wiUCP}, 
we obtain the following result generalizing \cite[Section 2, Corollary]{MP1}.

\begin{thm}\label{thm:cobD}
Let $X$ be an arbitrary topological space  
and $\D=(\D_x)_{\xX}$ be a family of Dirac type operators on a smooth compact Riemannian manifold $N$
	such that the coefficients of $\D_x$ and $\D^t_x$, 
	together with the first derivatives of the conormal symbol of $\D_x$ in the tangential directions, 
  depend continuously on $x$ in every local chart of $N$.
	Then the family $\AAt = (\AAt_x)$ of symmetrized boundary operators $\AAt_x = (\A_x+\A_x^t)/2$
	admits both a spectral section and a trivializing family, with smoothing trivializing operators.
\end{thm}

\begin{rem}\label{rem:even-dim}
While the statement of Theorem \ref{thm:cobD} does not distinguish between 
even- and odd-dimensional manifolds $N$,
the theorem is really interesting only in the even-dimensional case.

If a manifold is odd-dimensional, then its boundary is a closed manifold of even dimension.
At the same time, an \textit{arbitrary} family of symmetric Dirac type operators $\A=(\A_x)$ 
over a closed oriented manifold $M$ of dimension $2k$ admits a spectral section and a trivializing family,
regardless of whether $\A$ is cobordant to zero or not.

Indeed, let $A$ be a symmetric Dirac type operator over $M$ with the symbol $a$.
Let $\sigma$ be the \q{normalized orientation}, 
$\sigma_y = i^k\,a(\xi_1)\cdot\ldots\cdot a(\xi_{2k})$, 
where $(\xi_1,\ldots,\xi_{2k})$ is a positively oriented orthonormal basis of $T^*_y M$, $y\in M$.
Then $\sigma$ is a bundle symmetry anticommuting with $a$, so Theorem \ref{thm:ess-odd} can be applied.
In more concrete terms, the operators $\bar{A} = (A-\sigma A\sigma)/2$ and $A' = \bar{A}+\sigma$
are symmetric and have the same symbol as $A$.
Moreover, $\bar{A}$ anticommutes with $\sigma$,
so $(A')^2 = \bar{A}^2+1$ is invertible and thus $A'$ itself is invertible.
Therefore, every family $\A=(\A_x)$ of symmetric Dirac type operators over a closed even-dimensional manifold
is a bounded deformation of a family $\A' = (\A'_x)$ of invertible symmetric Dirac type operators.
The family $(P_x)$ of positive spectral projections $P_x = \chip(\A'_x)$
is a generalized spectral section for $\A$.
It follows from Theorems \ref{thm:ss} and \ref{thm:ss-tr-op} 
that $\A$ admits a spectral section and a trivializing family.

A relevant cobordism theorem in this case should take into account the grading $\sigma$ 
and state the existence of a $\Cl(1)$ spectral section for the family $\bar{\A}$ of odd operators.
We perform this in the next section, see Theorem \ref{thm:cobD1}.
\end{rem}

\section{Cobordism theorems: $\Z_2$-graded case}\label{sec:cob1} 

Let $N$ be a smooth compact Riemannian manifold with non-empty boundary
and $F$ be a smooth Hermitian vector bundle over $N$, as in the previous section.
Denote by $E$ the restriction of $F$ to the boundary $M = \pN$.


\sub{General cobordism theorem.}
Let $D$ be a first order \textit{symmetric} operator acting on sections of $F$, 
$A$ be its boundary operator, 
and $J = d(n)\in\Iso(E)$ be the conormal symbol of $D$.
Then $J$ is self-adjoint and the operator $JA+A^t J$ has zeroth order 
(that is, $JA+A^t J$ is a bundle endomorphism).
Suppose, in addition, that the \textit{symbol of $A$ is self-adjoint}. 
Then both $A-A^t$ and $JA+AJ$ are bundle endomorphisms.

Let $\sigma$ be a bundle symmetry (that is, a self-adjoint unitary bundle automorphism) of $E$ 
defined by the formula 
\begin{equation}\label{eq:JJ}
	\sigma = J\cdot|J|\inv = J\cdot(J^2)^{-1/2}.
\end{equation}
Then $\sigma A+A\sigma$ is also a bundle endomorphism. 
Indeed, the symbol $a$ of $A$ satisfies the anticommutation relation
\begin{equation}\label{eq:Jb}
	J(y) a(\xi)+a(\xi)J(y) = 0 \quad 
	 \text{for every } y\in M \text{ and } \xi\in T_y^*M.
\end{equation}
Thus the positive operator $T = J(y)^2$ commutes with $a(\xi)$.
Multiplying \eqref{eq:Jb} by $T^{-1/2}$,
we get $\sigma(y)a(\xi) + a(\xi)\sigma(y) = 0$,
which implies $\sigma A+A\sigma \in\End(E)$.

Instead of the symmetrized boundary operator $\At=(A+A^t)/2$,
we now consider the supersymmetrized operator
\begin{equation}\label{eq:Bsigma}
	\bar{A} = (\At - \sigma\At\sigma)/2 = (A + A^t - \sigma A\sigma - \sigma A^t\sigma)/4,
\end{equation}
which also differs from $A$ by a bundle endomorphism.
\index{supersymmetrized tangential operator}

In such a way, every first order symmetric elliptic operator $D$ acting on $F$
whose boundary operator has self-adjoint symbol
determines a $\Z_2$-grading $\sigma = \sigma_D$ of $E$
and an odd symmetric elliptic operator $\bar{A} = \bar{A}_D$ acting on $E$.

By \cite[Theorem 6.1.I]{BBLZ}, the space $\Lambda$ of Cauchy data of $D$ 
is a Lagrangian subspace of the symplectic Hilbert space $(H, iJ)$, $H = L^2(M;E)$.
In other words, the orthogonal projection $Q$ of $H$ onto $\Lambda$
(the Calder\'on projection)
satisfies the anticommutation relation $J(2Q-1) + (2Q-1)J = 0$.
Reasoning as above, we obtain $\sigma(2Q-1) + (2Q-1)\sigma = 0$.
Therefore, $Q$ is a generalized $\Cl(1)$ spectral section for $\bar{A}$,
as well as for any other odd symmetric operator with the symbol $a$.

Applying our previous results to this situation,
we obtain a graded version of the general cobordism theorem.

\begin{thm}\label{thm:cob1}
Let $X$ be an arbitrary topological space
and $\D=(\D_x)_{\xX}$ be a family of first order \textit{symmetric} elliptic differential operators 
acting on sections of $F$. 
Suppose that the principal symbol $a_x$ of the boundary operator of $\D_x$ is self-adjoint for every $\xX$
and that the conormal symbol $J_x$ of $\D_x$ and the Calder\'on projection $Q(\D_x)$ depend continuously on $x$.
Let the grading $\sigma_x$ on the boundary bundle $E$ be defined by the unitary part of $J_x$ as in \eqref{eq:JJ},
and let $\bar{\A} = (\bar{\A}_x)$ be a continuous family 
of odd symmetric differential operators $\bar{\A}_x$ over the boundary having $a_x$ as their principal symbols.
Then $\bar{\A}$ admits both a $\Cl(1)$ spectral section and a trivializing family of odd smoothing trivializing operators.
\end{thm}

\proof
The grading $\sigma$ determines a norm continuous map from $X$ to the space of symmetries
\[ \Sigma = \sett{S\in\Bsa}{S^2=1 \text{ and } \dim\ker(S+1) = \dim\ker(S-1) = \infty }, \] 
where $H=L^2(M;E)$.
The trivial Hilbert bundle over $\Sigma$ with the fiber $H$ 
has the canonical decomposition into the direct sum $\H^0\oplus\H^1$ 
of two Hilbert bundles, whose fibers over $S\in\Sigma$ are $\ker(S+1)$ and $\ker(S-1)$ respectively.
$\H^0$ and $\H^1$ are locally trivial bundles with infinite-dimensional fibers $H'$ and the structure group $\U(H')$.
The structure group 
is contractible by Kuiper's theorem \cite{Kui}.
The base space $\Sigma$ is metric and thus paracompact,
so both these bundles are trivial by \cite[Theorem 7.5]{Dold}.
The direct sum of their trivializations determines 
a norm continuous map $u\colon\Sigma\to\U(H)$ such that the symmetry $S_0 = u(S)\cdot S\cdot u(S)^*$ does not depend on $S$.

Conjugation by the map $v = u\circ\sigma\colon X\to\U(H)$ takes $\bar{\A}\colon X\to\Rsa_K(H)$
to a Riesz continuous map $\bar{\A}'\colon X\to\Rsa_K(H)$ 
and the family $x\mapsto Q(\D_x)$ of Calder\'on projections to a norm continuous map $Q'\colon X\to\Proj(H)$ 
so that all the operators $\bar{\A}'_x$ are odd and $Q'$ is a generalized $\Cl(1)$ spectral section for $\bar{\A}'$ 
with respect to the fixed grading $S_0$ of $H$.
Let $P'$ and $\C'$ be a $\Cl(1)$ spectral section and a trivializing family for $\bar{\A}'$ provided by Theorem \ref{thm:ss1}.
Then $P_x=v_x^*P'_x v_x$ is a $\Cl(1)$ spectral section 
and $\C_x=v_x^*\C'_x v_x$ is a trivializing family for $\bar{\A}$
with respect to the family $(\sigma_x)$ of gradings.
Since the range of $\C_x$ lies in the span of a finite number of $C^{\infty}$-sections 
(namely, eigenfunctions of $\bar{\A}_x$), all the operators $\C_x$ are smoothing. 
\endproof

\sub{Dirac type operators.}
If $D$ is a symmetric Dirac type operator, then the conormal symbol of $D$ is unitary 
and thus coincides with the grading $\sigma$ defined by \eqref{eq:JJ}.
Applying Theorem \ref{thm:cob1} to this situation,
we obtain the following generalization of \cite[Corollary 1]{MP2}. 

\begin{thm}\label{thm:cobD1}
  Let $X$ be an arbitrary topological space and 
  $\D=(\D_x)_{\xX}$ be a family of \textit{symmetric} Dirac type operators satisfying conditions of Theorem \ref{thm:cobD}.	
	Then the family $\bar{\A} = (\bar{\A}_x)$ of supersymmetrized boundary operators 
	admits both a $\Cl(1)$ spectral section and 
	a trivializing family with odd smoothing trivializing operators,
	with respect to the family of gradings of boundary bundles given by the conormal symbol of $\D_x$.
\end{thm}

\proof
By Propositions \ref{prop:Cald-cont} and \ref{prop:wiUCP},
the Calder\'on projection $Q(\D_x)$ depends continuously on $x$.
By the definition of the supersymmetrized boundary operator, 
$\bar{\A}_x$ depends continuously on $x$ 
and its symbol coincides with the symbol $a_x$ of the boundary operator of $\D_x$.
It remains to apply Theorem \ref{thm:cob1}.
\endproof


\addtocontents{toc}{\vspace{-5pt}}

\part{Graph continuous families}\label{part:graph}

\addtocontents{toc}{\vspace{-12pt}}

Throughout this part, all families of regular operators are supposed to be graph continuous.

Recall that the Cayley transform of a regular self-adjoint operator $A$ 
is the unitary operator defined by the formula 
$\kappa(A) = (A-i)(A+i)\inv$.
\index{2@Maps!2k1@$\kappa$, Cayley transform}
\index{Cayley transform}
The Cayley transform $\kappa\colon{^g\Reg\sa}(H)\to\U(H)$ is a homeomorphism on the image \cite[Theorem 1.1]{BBLP}.
Therefore, the graph topology on the subspace $\Reg\sa(H)$ of $\Reg(H)$ 
can be equivalently described as the topology induced by the inclusion 
$\kappa\colon\Reg\sa(H)\hookto\U(H)$.
We will use this fact below.

\section{Semibounded operators}\label{sec:gt-sb}

\upskip

\sub{Positive operators.}
Let $\Reg^+(H)$ denote the subspace of $\Reg\sa(H)$ consisting of positive operators.
\index{1@Spaces!1R3+@$\Reg^+(H)$, positive regular operators}

\begin{prop}\label{prop:bdd}
  The restrictions of the graph topology and the Riesz topology to $\Reg^+(H)$ coincide.
\end{prop}

\proof
For $a=f(A)$ the identity $(1+A^2)\inv = 1-a^2$ implies
\[ \kappa(A) = \frac{A-i}{A+i} = \frac{a-i\sqrt{1-a^2}}{a+i\sqrt{1-a^2}} = \br{a-i\sqrt{1-a^2}}^2 = \kat(a), \]
where $\kat\colon[-1,1]\to\U(\CC) = \sett{z\in\CC}{|z|=1}$ is a continuous function 
given by the formula $\kat(a) = \br{a-i\sqrt{1-a^2}}^2$.
\index{2@Maps!2k2@$\kat$}
Thus the Cayley transform factors through the bounded transform: 
$\kappa = \kat\circ f$.
The function $\kat$ is not invertible.
However, its restriction to the interval $[0,1]$ \textit{is} invertible: 
it is a homeomorphism from $[0,1]$ to the bottom half of the unite circle 
$\Gamma = \sett{e^{it}}{t\in [-\pi,0] }$.
The inverse homeomorphism $\varphi\colon\Gamma\to[0,1]$ is given by the formula $\varphi(e^{it})=\cos(t/2)$. 
Hence the restriction of the bounded transform $f$ to $\Reg^+(H)$
coincides with the composition $\varphi\circ\kappa$ and thus is continuous.
It follows that the restriction of the graph topology to $\Reg^+(H)$ coincides with the Riesz topology.
\endproof

\sub{Semibounded operators.}
An operator $A\in\Rsa(H)$ is called bounded from below if $A-c$ is positive for some $c\in\R$.
Similarly, $A$ is called bounded from above if $c-A$ is positive for some $c\in\R$
(equivalently, $-A$ is bounded from below).
$A$ is called semibounded if it is bounded either from below or from above.

Proposition \ref{prop:bdd} can be easily generalized to operators 
bounded from below or from above by a \textit{fixed} constant $c$.

\begin{prop}\label{prop:sbd}
The graph topology coincides with the Riesz topology on the subsets
\[
 \Reg^{\geq c}(H) = \sett{A\in\Reg\sa(H)}{A-c\geq 0 } \quad\text{and}\quad
 \Reg^{\leq c}(H) = \sett{A\in\Reg\sa(H)}{c-A\geq 0}
\]
of $\Reg\sa(H)$ for every $c\in\R$.

\end{prop}

\proof
We work as in the proof of Proposition \ref{prop:bdd}, with the same designations.
The function $\kat$ defines a homeomorphism from the interval $[f(c),1]$ to the arc 
$\Gamma_c = \sett{e^{it}}{t\in [t_c,0] }$, where $t_c\in(-2\pi,0)$, $e^{it_c}=\kappa(c)$.
The inverse homeomorphism $\varphi\colon\Gamma_c\to[f(c),1]$ is given by the same formula as before,
$\varphi(e^{it})=\cos(t/2)$. 
The restriction of the bounded transform $f$ to $\Reg^{\geq c}(H)$
coincides with the composition $\varphi\circ\kappa$ and thus is continuous.
It follows that the restriction of the graph topology to $\Reg^{\geq c}(H)$ coincides with the Riesz topology.

The proof for $\Reg^{\leq c}(H)$ is quite similar 
(or, equivalently, can be obtain from the result for $\Reg^{\geq -c}(H)$ using the map $A\mapsto -A$,
which is a homeomorphism both in the graph and the Riesz topology).
\endproof

\sub{Families of semibounded operators.}
By \cite[Addendum, Theorem 1]{CL}, the restriction of the graph topology
to the subspace of bounded operators coincides with the usual norm topology, 
and thus with the Riesz topology.

This property cannot be generalized to the subspace of \textit{all semibounded} 
self-adjoint regular operators. 
Indeed, Examples \ref{ex:delta} and \ref{ex:Hl2N} below demonstrate families of semibounded operators 
(namely, each operator is bounded from below), 
which are graph continuous but Riesz discontinuous. 
Riesz discontinuity there is caused by the absence of a \textit{continuous} lower bound,
as the following proposition shows.

\begin{prop}\label{prop:sbd2}
Let $\A\colon X\to\Rsa_K(H)$ be a graph continuous family of regular self-adjoint operators, 
with each $\A_x$ bounded from below.
Then the following two conditions are equivalent:
\begin{enumerate}[topsep=-3pt, itemsep=0pt, parsep=3pt, partopsep=0pt]
	\item $\A$ is Riesz continuous.
	\item There is a continuous function $c\colon X\to\R$
             such that $\A_x-c_x\geq 0$ for every $\xX$.
\end{enumerate}
The analogous equivalence holds for operators bounded from above.
\end{prop}

\proof
($1\Rightarrow 2$)
Let $c_x\in\R$ be the exact lower bound of $\A_x$. 
Then $\bar{c}_x = f(c_x)\in(-1,1)$ is the exact lower bound of $a_x=f(\A_x)$.
If $\A$ is Riesz continuous, then $x\mapsto a_x$ is norm continuous,
and thus $\bar{c} = f\circ c$ is a continuous function from $X$ to $(-1,1)$.
Therefore, $c = f\inv\circ \bar{c} \colon X\to\R$, 
$c_x = \bar{c}_x\br{1-\bar{c}_x^2}^{-1/2}$, is also continuous.

($2\Rightarrow 1$)
Let $c$ be such a function.
Then every point $\xX$ has an open neighborhood $U_x$ such that $c_y > c_x-1$ for every $y\in U_x$.
The restriction of $\A$ to $U_x$ satisfies conditions of Proposition \ref{prop:sbd}
(with the constant lower bound $c_x-1$) and thus is Riesz continuous.
Since the open subsets $U_x$ cover $X$, $\A$ itself is Riesz continuous.

The same reasoning but with inverted signs works for operators bounded from above.
This completes the proof of the proposition.
\endproof

\section{Spectral sections}\label{sec:gt-ss}

\upskip
\sub{Arbitrary base spaces.}
It follows from Theorem \ref{thm:ss} that
a Riesz continuous family always has local spectral sections locally,
and the only obstruction for existence of a global spectral section is a topological one. 
In contrast with this, 
a graph continuous family may have no spectral section even locally.
In fact, Riesz continuity is \textit{necessary} for a local existence of a spectral section,
as the following result shows.

\begin{thm}\label{thm:ss2Riesz}
  Let $X$ be an arbitrary topological space and 
	$\A\colon X\to\CRR\sa(H)$ be a graph continuous map having a spectral section.
	Then $\A$ is Riesz continuous.
\end{thm}

\proof
Let $P$ be a spectral section for $\A$ and $r\colon X\to\R_+$ be the corresponding cut-off function.
Let $x_0\in X$.
Choose a positive constant $r_0$ and a neighborhood $V$ of $x_0$ such that 
$\pm r_0\in\Res(\A_{x})$ and $r_0>r(x)$ for all $x\in V$.

The finite rank projection $\So_x = \one_{(-r_0,r_0)}(\A_x)$ 
depends continuously on $x\in V$ and commutes with $P_x$.
Hence 
\begin{equation}\label{eq:ss2Riesz}
	\So_x, \quad S^+_x = (1-\So_x)P_x, \;\text{ and }\; S^-_x = (1-\So_x)(1-P_x)
\end{equation}
are mutually orthogonal projections depending continuously on $x\in V$.
Decreasing $V$ if necessary, one can find a continuous map $g\colon V\to\U(H)$ 
such that the conjugation by $g$ takes these three projection-valued maps 
to constant projections $\So$, $S^+$, and $S^-$.
Indeed, one can first find a neighborhood $V_1\subset V$ of $x_0$
and a map $g_1\colon V_1\to\U(H)$ such that $g_1(x) \So_x g_1(x)^* \equiv \So$ 
\cite[Proposition 5.2.6]{WO}.
Next, one can find a neighbourhood $V_2\subset V_1$ of $x_0$ 
and a map $g_2\colon V_2\to\U(H)$ such that 
$g_2(x)$ is equal to the identity on the range of $\So$
and the conjugation by $g_2(x)$ takes $g_1(x)S^+_x g_1(x)^*$ to $S^+$.
Then $g = g_2 g_1$ is a desired trivialization of projections \eqref{eq:ss2Riesz} over $V_2$.

Let $H^{\circ}$, $H^+$, and $H^-$ be the ranges of projections $\So$, $S^+$, and $S^-$. 
Then
\begin{equation}\label{eq:AAA}
	g_x\,\A_x\,g_x^* = \A^-_x\oplus \A^{\circ}_x\oplus \A^+_x
\end{equation}
with respect to the orthogonal decomposition $H = H^-\oplus H^{\circ} \oplus H^+$.
The map $g\A g^*\colon V\to\CRR\sa(H)$ is graph continuous, 
so its components $\A^-$, $\A^{\circ}$, $\A^+$ are also graph continuous.
Since $\A^{\circ}$ acts on the finite-dimensional space $\Ho$, it is norm continuous.
Both $\A^+$ and $-\A^-$ are bounded from below by a positive constant $r_0$, 
so by Proposition \ref{prop:bdd} they are Riesz continuous.
Substituting this to \eqref{eq:AAA}, we see that $g\A g^*$ is Riesz continuous on $V$,
and thus the restriction of $\A$ to $V$ is also Riesz continuous.
Since $x_0\in X$ was chosen arbitrarily, $\A$ is Riesz continuous on the whole $X$.
This completes the proof of the theorem.
\endproof

\sub{Compact base spaces.}
Taking together Theorem \ref{thm:ss2Riesz} and Proposition \ref{prop:MP1},
we immediately obtain the following result.

\begin{thm}\label{thm:comp-gap}
  Let $X$ be a compact space and $\A\colon X\to\CRR\sa(H)$ be a graph continuous map.
	Then the following two conditions are equivalent:
\begin{enumerate}[topsep=-3pt, itemsep=0pt, parsep=3pt, partopsep=0pt]
	\item $\A$ has a spectral section.
	\item $\A$ is Riesz continuous and $\ind(\A) = 0 \in K^1(X)$.
\end{enumerate}
\end{thm}

\begin{rem}
Theorems \ref{thm:ss2Riesz} and \ref{thm:comp-gap} show that the straightforward transfer of 
\cite[Proposition 1]{MP1} and \cite[Proposition 2]{MP2}
to elliptic operators on manifolds with boundary does not work; 
one needs to be very careful using spectral sections in this framework.
For example, J.~Yu applies \cite[Proposition 1]{MP1} 
to families of Dirac operators with local boundary conditions in \cite{Yu}.
However, to make such an application justified,
one needs to ensure first that the corresponding families of unbounded operators are Riesz continuous.
This is typically a very non-trivial problem, 
the answer to which seems to be unknown in the general situation considered in \cite{Yu}. 
It is unknown, except for several special cases (see e.g. \cite{BR}), 
whether families of unbounded operators defined by boundary value problems are Riesz continuous.

The following example of Rellich \cite[Example V-4.14]{Kato} 
demonstrates lack of Riesz continuity already in the one-dimensional framework.
\end{rem}

\begin{exm}[Graph continuous but Riesz discontinuous family of local boundary value problems]\label{ex:delta}
Let $A = -d^2/d_t^2$ be the second order operator acting on complex-valued functions 
$\psi\colon M=[0,1]\to\CC$.
For $x\in\R$, let $\A_x$ be the operator $A$ with the domain 
\begin{equation}\label{eq:delta}
	\dom(\A_x) = \sett{\psi\in H^2(M;\CC)}{\psi(0) = 0, \; \psi(1) = x\psi'(1)}.
\end{equation}
given by local boundary conditions.
Here the type of the boundary condition on the right end of the interval $[0,1]$
changes from Dirichlet boundary condition at $x=0$ to mixed boundary condition at $x\neq 0$.

Every $\A_x$ considered as an operator on $H = L^2(M,\CC)$ is a self-adjoint operator with compact resolvent.
The map $H^2(M;\CC)\to\CC^3$, $\psi\mapsto(\psi(0),\psi(1),\psi'(1))$, is continuous,
so the domain of $\A_x$ is a closed subspace of $H^2(M;\CC)$
depending continuously on $x$ in the gap topology on $\Gr\br{H^2(M;\CC)}$.
It follows from \cite[Proposition A.9]{Pr17} that
the family $\A=(\A_x)$ of unbounded operators on $H$ is graph continuous.

Each operator $\A_x$ is bounded from below.
It can be easily seen that $\A_x$ has a negative eigenvalue if and only if $x\in(0,1)$.
The function $\psi(t) = e^{\mu t} - e^{-\mu t}$, $\mu>0$, is an eigenfunction of $\A_x$
if $e^{2\mu}-1 = x\mu(e^{2\mu}+1)$;
the corresponding eigenvalue is $\lambda = -\mu^2$.
As $x$ goes to $+0$, $\mu$ goes to $+\infty$ and thus $\lambda$ goes to $-\infty$.
It follows that the negative eigenvalue of $a_x=f(\A_x)$ goes to $-1$ as $x\to +0$.
However, $A_0$ and thus $a_0$ is positive.
Therefore, $a=f\circ\A$ is discontinuous at $x=0$ and thus $\A$ is Riesz discontinuous at $x=0$.

In conclusion, the map $\A\colon\R\to\CRR\sa(H)$ is graph continuous,
but its restriction to every interval $[0,\eps]$ is Riesz discontinuous
and does not admit a spectral section.
\end{exm}

\section{Generalized spectral sections}\label{sec:gt-gss}

In the previous section we showed that a graph continuous family 
admitting a spectral section has to be Riesz continuous.
Situation with generalized spectral sections, however, is more ambiguous. 
We give here several illustrating examples.

\sub{Graph continuous but Riesz discontinuous family with generalized spectral section.}
In Example \ref{ex:delta} above a graph continuous family of local boundary value problems
is Riesz discontinuous but has a generalized spectral section $P\equiv 1$.
The following example of Fuglede (presented in \cite[Remark 1.6]{Nic07}) demonstrates a discrete analogue of such situation.

\begin{exm}\label{ex:Hl2N}
Let $X=\N\cup\set{\infty}$ be a one-point compactification of $\N$.
Let $(e_n)_{n\in\N}$ be an orthonormal basis of $H$.
Consider a family $\A\colon X\to\CRRsa$ of diagonal (in the chosen basis) 
semibounded operators given by the formulae
\[ \A_x(e_n) = \case{n, & n\neq x \\ -n, & n=x} \;\text{ for } x\in\N; \quad \A_{\infty}(e_n) = n. \]
Since $\norm{\kappa(\A_{\infty}) - \kappa(\A_x)} = \left|\kappa(x) - \kappa(-x)\right| \to 0$ 
as $x\to\infty$,
the family $\A$ is graph continuous.
On the other side, $\norm{f(\A_{\infty}) - f(A_x)} = 2f(x) \to 2$ as $x\to\infty$,
so $\A$ is Riesz discontinuous at $\infty$.

The constant function $P\colon X\to\Proj(H)$ taking every $\xX$ to the identity 
is a generalized spectral section for $\A$.
Moreover, $P_x$ is even an $r_x$-spectral section for $\A_x$, 
where $r\colon X\to\R_+$ is an arbitrary function such that $r_x > x$ for $x\in\N$.
However, every such function $r$ is discontinuous at $\infty$,
so $P$ is not a global spectral section for $\A$. 
(Otherwise, of course, we would have a contradiction with Theorem \ref{thm:ss2Riesz}.)
\end{exm}

\upskip
\sub{Graph continuous families without generalized spectral sections.}
A generalized spectral section does not necessarily exist
even for a contractible base space (Example \ref{ex:pm})
or for a family of invertible semibounded operators (Example \ref{ex:Hl2N2}).

\begin{exm}\label{ex:pm}
The space $\CRRsa$ equipped with the graph topology is path connected \cite{Jo}.
Let $\A\colon[0,1]\to\CRRsa$ be a graph continuous path connecting 
a negative operator $\A_0$ with a positive operator $\A_1$.
(Such a path $\A$ can even be chosen consisting of invertible operators, but we do not explore it here for simplicity.)
Then $\A$ has no generalized spectral section.	 	
Indeed, a generalized spectral section $P_0$ for $\A_0$ should be compact,
while a generalized spectral section $P_1$ for $\A_1$ should have compact complement $1-P_1$.
Any two such projections $P_0$ and $P_1$ lie in the different connected components of the space $\Proj(H)$,
so they cannot be connected by a path $P\colon [0,1]\to\Proj(H)$.
\end{exm}

\begin{exm}\label{ex:Hl2N2}
Let $X=\N\cup\set{\infty}$ be a one-point compactification of $\N$
and $(e_n)_{n\in\N}$ be an orthonormal basis of $H$, as in example \ref{ex:Hl2N}.
Consider a family $\A\colon X\to\CRRsa$ of invertible semibounded operators given by the formulae
\[ \A_x(e_n) = \case{-n, & n<x \\ \;\; n, & n\geq x} \;\text{ for } x\in\N; \quad \A_{\infty}(e_n) = -n.\]
The same reasoning as in example \ref{ex:Hl2N} shows that $\A$ is graph continuous. 
However, $\A$ has no generalized spectral section. 
Indeed, the operator $\A_{\infty}$ is bounded from above,
while all the other $\A_x$, $x\in\N$, are bounded from below.
Thus a generalized spectral section $P_{\infty}$ for $\A_{\infty}$ should be compact,
while a generalized spectral section $P_x$ for $\A_x$, $x\in\N$, 
should have compact complement $1-P_x$.
Since $P_{\infty}$ cannot be a limit point of $P_x$, 
such a map $P$ cannot be norm continuous.
\end{exm}

{\small\input{Spectral-sections.ind}}

\end{document}